\documentclass[MSNbibl,number,citesort,rotating,seceqn,dvips]{arxbj}
\usepackage{graphicx}

\aid{0}
\volume{19}
\issue{5B}
\pubyear{2013}
\firstpage{2590}
\lastpage{2626}
\doi{10.3150/12-BEJ467} 

\makeatletter
\define@key{bid}{pmid}{}

\newremark{remark}{Remark}[section]
\newtheorem{theorem}{Theorem}[section]
\newproclaim{definition}{Definition}[section]
\newtheorem{proposition}{Proposition}[section]

\newcommand{\overset}{\stackrel}
\newcommand{\Gam}{\operatorname{Gamma}}
\newcommand{\Bet}{\operatorname{Beta}}
\renewcommand{\emptyset}{\varnothing}
\makeatother

\begin{document}
\begin{frontmatter}

\title{A conjugate class of random probability measures based on
tilting and with its posterior analysis}
\runtitle{A conjugate class of random probability measures}

\begin{aug}
\author{\fnms{John W.} \snm{Lau}\corref{}\ead[label=e1]{john.lau@uwa.edu.au}}
\runauthor{J.W. Lau} 
\address{School of Mathematics and Statistics, The University of Western
Australia, 35 Stirling Highway, Crawley, WA 6009, Australia. \\\printead{e1}}
\end{aug}

\received{\smonth{10} \syear{2010}}
\revised{\smonth{5} \syear{2012}}

%
\begin{abstract}
This article constructs a class of random probability measures based
on exponentially and polynomially tilting operated on the laws of
completely random measures. The class is proved to be conjugate in
that it covers both prior and posterior random probability measures
in the Bayesian sense. Moreover, the class includes some common
and widely used random probability measures, the normalized
completely random measures (James (Poisson process partition calculus
with applications to exchangeable models and Bayesian nonparametrics
(2002) Preprint),
Regazzini, Lijoi and Pr\"unster (\textit{Ann. Statist.} \textbf{31}
(2003) 560--585),
Lijoi, Mena and Pr\"unster (\textit{J.~Amer. Statist. Assoc.} \textbf
{100} (2005) 1278--1291)) and the Poisson--Dirichlet process
(Pitman and Yor (\textit{Ann. Probab.} \textbf{25} (1997) 855--900),
Ishwaran and James (\textit{J. Amer. Statist. Assoc.} \textbf{96}
(2001) 161--173),
Pitman (In \textit{Science and Statistics: A Festschrift for Terry
Speed} (2003) 1--34 IMS)), in a
single construction. We describe an augmented version of the
Blackwell--MacQueen P\'olya urn sampling scheme (Blackwell and MacQueen
(\textit{Ann. Statist.} \textbf{1} (1973) 353--355))
that simplifies implementation and provide a simulation study for
approximating the probabilities of partition sizes.
\end{abstract}

%

\end{frontmatter}

\section{Introduction}

Random probability measures derived from normalized independent
increment processes have been studied for decades.
Kingman \cite{Kingman1975} considers normalization over subordinators of
L\'evy processes with only positive jumps. Regazzini \textit{et al.} \cite{RLP2003} introduce
the class of the normalized random measures with independent
increments (normalized {\it completely random measures}) for
studying the probabilistic properties of mean functionals of random
probability measures. Investigations for statistical modeling are
available in James \cite{James2002}, Lijoi \textit{et al.} \cite{LMP2005b} and
James \textit{et al.} \cite{JLP2008}.

In Bayesian non-parametric statistics, the normalized random process
is considered to be an unknown parameter and the posterior
distribution of the process is usually of interest. The most popular
class of such random processes for statistical modelling is the
Dirichlet process (Ferguson \cite{Ferguson1973}; Lo \cite{Lo1984}). The
Dirichlet process is appealing because it induces model flexibility
and it is also conjugate in the sense that the posterior process,
the process conditional on the data, is also a Dirichlet process. In
fact, a surprising result shown in James \textit{et al.} \cite{JLP2006} states that only
the Dirichlet process has the conjugacy property among the random
probability measures in the class of normalized independent
increment processes. In the present work, we are able to show that
it is not the case for a richer class of normalized processes. The
class of normalized processes considered in this article is derived
based on tilting, in particular, exponentially and polynomially
tilting operated
on the laws of completely random measures
and this yields the class of laws containing the prior random
probability measures and their posteriors in the Bayesian sense.
So, the random probabilities in this class is conjugate.

Tilting the laws of random processes provides a way to enrich the
class of random processes through change of measure. For example,
Pitman \cite{Pitman2003} constructs the Poisson--Kingman process by
normalizing a random process that has a tilted law of completely
random measures. Other well known special cases are the
Poisson--Dirichlet process, whose law is constructed by polynomially
tilting the laws of positive $\alpha$-stable processes
(Pitman \cite{Pitman2003}), and the beta-gamma process
(James \cite{James2005b}), whose law is obtained by polynomially
tilting the laws of gamma processes. However, these studies give
less attention to the statistical properties of the tilted
processes. Our studies in fact focus on showing conjugate property
of the random probability measure derived from tilted laws of the
completely random measures and providing the posterior analysis of
the class of random probability measures.

Applications of non-parametric models are becoming increasingly
common in Bayesian statistics. However, implementing non-parametric
models is rarely straightforward and often involves Markov chain
Monte Carlo (MCMC) algorithms that might require evaluation of
complicated functions. 
This article provides an augmented form of the sampling algorithm,
namely the Blackwell--MacQueen P\'olya urn sampling scheme, which
avoids the necessity of evaluating those functions for some special
cases and which we believe would be beneficial to the use of
normalized random measures in statistical applications in future. We
provide a simulation study concerned with estimating probabilities
of partition sizes,
a problem that arises in biological speciation (Lijoi \textit{et al.} \cite{LMP2007a}). We
show our algorithm
yields similar results to three other sampling schemes based on the
Blackwell--MacQueen P\'olya urn distribution.

The article proceeds as follows. Section \ref{SectionConstruction}
describes the construction of a class of random probability measures
obtained by tilting the laws of completely random measures.
Materials presented here are in compact form; we refer the reader to
Daley and Vere-Jones \cite{DV2008} for the complete treatment on these topics (see also
Kingman \cite{Kingman1967,Kingman1975,Kingman1993};
Kallenberg \cite{Kallenberg1983,Kallenberg2001}).
Section \ref{SectionNRM-PET} considers a class of random
probability measures constructed through polynomially and
exponentially tilting the laws of completely random measures.
Section \ref{SectionNRM-PET} also provides details on the prior
and posterior distributions, proves the conjugacy property and
describes the augmented Blackwell--MacQueen sampling scheme. Section
\ref{Sectiontiltedexamples} describes two specific examples of the
tilting strategy, tilting the laws of the generalized gamma process and
the generalized
Dirichlet process is demonstrated. Section \ref{Sectionsim} presents
the simulation
study. Section \ref{SectionConclusionRemarks} concludes the
articles and provides a future research perspective. Proofs of
theorems and propositions are included in the \hyperref[app]{Appendix}.

In addition, as a referee pointed out that theorems in Section
\ref{SectionNRM-PET} are not entirely new and a version of them
has appeared in an unpublished manuscript (James \cite{James2002}), we
would like to acknowledge that the class of random probability
measures involving polynomially and exponentially tilting was first
studied in the manuscript (James \cite{James2002}, Chapter 5),
including a posterior analysis. We will further comments on these
aspects and connections in Section \ref{SectionNRM-PET}.

\section{Construction of the class of random probability measures}\label{SectionConstruction}%

Let the triple $(\Omega,\mathfrak{F},\textsf{P})$ be the basic
probability space. Assume $\mathbb{X}$ is a Polish space endowed
with a metric $d_{\mathbb{X}}$ generating a Borel $\sigma$-field
$\mathcal{B}(\mathbb{X})$.
Let $\mathcal{M}_\mathbb{X}$
denote the space of boundedly finite measures on
$(\mathbb{X},\mathcal{B}(\mathbb{X}))$. A measure is said to be
boundedly finite if it is finite on bounded sets. The space
$\mathcal{M}_\mathbb{X}$ is a Polish space equipped with the metric
of weak convergence. This induces the Borel $\sigma$-field
$\mathcal{B}(\mathcal{M}_\mathbb{X})$.
A {\it random measure}, $\mu$ say, taking values from
$\mathcal{M}_\mathbb{X}$, is a measureable mapping
$\mu\dvtx \Omega\times\mathbb{X}\rightarrow\mathbb{R}^+\cup\{0\}$ where
$\mathbb{R}^+$ denotes the positive real line. For each
$\omega\in\Omega$, $\mu(\omega,\cdot)$ is a boundedly finite measure
on $(\mathbb{X},\mathcal{B}(\mathbb{X}))$ and $\mu(\omega,A)$ is a
random variable for all bounded sets $A\in\mathcal{B}(\mathbb{X})$.
For convenience of notation, we write $\mu(A)$ instead of
$\mu(\omega,A)$ from now on. For further details, see
Daley and Vere-Jones \cite{DV2008}, Chapter 9.


A random measure $\mu$ is a {\it completely random
measure} ({\sc crm}) on the measure space
$(\mathbb{X},\mathcal{B}(\mathbb{X}))$ if for all finite families of
pairwise disjoint, bounded Borel sets $A_1, A_2,\ldots, A_k
\in \mathcal{B}(\mathbb{X})$, the random variables $\mu(A_1)$,
$\mu(A_2), \ldots, \mu(A_k)$ are mutually independent.
For any {\sc crm}, there is a representation theorem due to
Kingman \cite{Kingman1967}, Theorem 1 (see also Kingman \cite{Kingman1975} and
Kingman \cite{Kingman1993}, Chapter 8), and the theorem is nicely
described in Daley and Vere-Jones \cite{DV2008}, Theorem~10.1.III. The version in Daley and Vere-Jones \cite{DV2008}, Theorem 10.1.III, says that
a {\sc crm} $\mu$ can be represented as a sum of an atomic measure
with countably many fixed atoms, a deterministic non-atomic measure
and a measure derived from a Poisson process. The representation is
given by
%
\begin{eqnarray}
\label{representation} \mu(A)=\sum_{k=1}^{\infty}
U_k\delta_{x_k}(A)+\lambda(A)+\int_{\mathbb{R}^+}
y N(A,\mathrm{d}y),\qquad A\in\mathcal{B}(\mathbb{X}),
\end{eqnarray}
where the sequence $\{x_1, x_2, \ldots \}$ is the countable
set of fixed atoms of $\mu$, $\{U_1, U_2, \ldots \}$ is a
sequence of mutually independent non-negative random variables,
$\lambda$ is a fixed non-atomic boundedly finite measure on
$\mathbb{X}$, and $N$ is a Poisson process on
$\mathbb{X}\times\mathbb{R}^+$. This Poisson process $N$ is independent
of $\{U_1, U_2,
\ldots \}$ and has an intensity measure $\nu$ on $\mathbb{X}\times
\mathbb{R}^+$. The intensity
$\nu$ satisfies the following two conditions. For every bounded set~$A\in\mathcal{B}(\mathbb{X})$,
%
\begin{eqnarray}
\label{condition1} \int_{\mathbb{R}^{+}} \min\{y,s\} \nu(A,\mathrm{d}y) <
\infty,\qquad s\in\mathbb{R}^+,
\end{eqnarray}
and
%
\begin{eqnarray}
\label{condition2} \nu\bigl(\{x\},\mathbb{R}^+\bigr)=0,\qquad x\in\mathbb{X}.
\end{eqnarray}
%
Notice that conditions (\ref{condition1}) and (\ref{condition2}) guarantee
that the random measure $\int_{\mathbb{R}^+} y N(\cdot,\mathrm{d}y)$ is
boundedly finite on $\mathbb{X}$ and has no fixed atoms respectively
(see also
Kallenberg \cite{Kallenberg2001}, Chapter 12).

The Poisson process in (\ref{representation}) can be considered as a
marked Poisson process on $\mathbb{X}$ having mark space
$\mathbb{R}^+$.
Again the product space $\mathbb{X}\times\mathbb{R}^+$ is a Polish
space with a suitable metric $d_{\mathbb{X}\times\mathbb{R}^+}$
extended from $d_{\mathbb{X}}$. The Borel $\sigma$-field of the
Polish space $\mathbb{X}\times\mathbb{R}^+$ is given by
$\mathcal{B}(\mathbb{X}\times\mathbb{R}^+)=\mathcal{B}(\mathbb
{X})\otimes\mathcal{B}(\mathbb{R}^+)$
where $\mathcal{B}(\mathbb{R}^+)$ denotes the $\sigma$-field
generated by the open subsets of $\mathbb{R}^+$. Then the Poisson
process $N$ is a mapping
$N\dvtx \Omega\times\mathbb{X}\times\mathbb{R}^+\rightarrow
\mathbb{Z}^+\cup\{0\}$ where $\mathbb{Z}^+$ denotes the positive
integers. This Poisson process takes value from the space of
boundedly finite measures
$(\mathcal{M}_{\mathbb{X}\times\mathbb{R}^+},\mathcal{B}(\mathcal
{M}_{\mathbb{X}\times\mathbb{R}^+}))$
defined analogous to the space
$(\mathcal{M}_{\mathbb{X}},\mathcal{B}(\mathcal{M}_{\mathbb{X}}))$
discussed in the first paragraph. Here we write $N(A,B)$ to
represent $N(\omega,A,B)$ for convenience of notation. The intensity
measure $\nu$ of this Poisson process is a non-atomic
$\sigma$-finite measure $\nu\dvtx \mathbb{X}\times\mathbb{R}^+\rightarrow
\mathbb{R}^+\cup\{0\}$.
The intensity
measure is particularly important since it is the only parameter of
the random measure and also the first moment of the measure. It
determines the nature of the process and further it also determines
the nature of the random measures derived from the Poisson process.
Here (\ref{condition2}) also implies that $N(\cdot,\mathbb{R}^+)$ is a
simple point process on $\mathbb{X}$. This point process
$N(\cdot,\mathbb{R}^+)$ is called ground process in Daley and Vere-Jones \cite{DV2008}, Chapter 9.

Initially, we restrict our attention to the {\sc crm} on
$\mathbb{X}$ without the first two components, the atomic component
and the drift, in (\ref{representation}), that is the {\sc crm} is
in the form of
%
\begin{eqnarray}
\label{CRM} \widetilde{\mu}(A):=\int_{\mathbb{R}^+} y N(A,\mathrm{d}y),\qquad A\in
\mathcal{B}(\mathbb{X}),
\end{eqnarray}
with respect to the Poisson process $N$ defined on
$\mathbb{X}\times\mathbb{R}^+$ with intensity measure $\nu$
satisfying conditions (\ref{condition1}) and (\ref{condition2}). The
law of the {\sc crm} $\widetilde{\mu}$, denoted by
$\mathcal{P}_{\widetilde{\mu}}$, which is derived from the law of
the Poisson process. For the sake of simplicity, we say
$\widetilde{\mu}$ has the parameter measure $\nu$. The measure $\nu$
can be decomposed into two measures, $\rho_x$ and $\eta$, written as
$\nu(\mathrm{d}x,\mathrm{d}s)=\eta(\mathrm{d}x)\rho_x(\mathrm{d}s)$. Such a decomposition is guaranteed
by Kallenberg \cite{Kallenberg1983}, Appendix~15.3.3, in which the measure $\rho_x$ is
uniquely determined outside any set of $\nu$ measure zero. Here
$\rho_x$ is a mapping
$\rho_x\dvtx \mathbb{R}^+\rightarrow\mathbb{R}^+\cup\{0\}$ for any
$x\in\mathbb{X}$ such that $\rho_x(A)$ is $\mathbb{X}$ measureable
for every bounded set $A\in\mathcal{B}(\mathbb{R}^+)$ and $\rho_x$
is a $\sigma$-finite measure. In particular, when $\rho_x$ is
dependent on $x\in\mathbb{X}$, the {\sc crm} $\widetilde\mu$ is
non-homogeneous. Otherwise, when $\rho_x$ is not dependent on
$x\in\mathbb{X}$, the {\sc crm} $\widetilde\mu$ is homogeneous. The
$\sigma$-finiteness of $\rho_x$ ensures the {\sc crm}
$\widetilde{\mu}$ has countably infinite jumps on any bounded set in
$\mathcal{B}(\mathbb{X})$. Here $\eta$ is a finite non-atomic
measure $\eta\dvtx \mathbb{X}\rightarrow\mathbb{R}^+\cup\{0\}$. Without
loss of generality, the measure $\eta$ is restricted to be a proper
probability measure on $\mathbb{X}$. This implies that the total
mass of the measure $\widetilde{\mu}$ is finite almost surely. Then
a random probability measure could be defined according to the ratio
of $\widetilde{\mu}(A)$ and the total mass $\widetilde{\mu}(\mathbb
{X})$, that is
$\widetilde{G}(A):=\widetilde{\mu}(A) /\widetilde{\mu}(\mathbb{X})$
for $A\in\mathcal{B}(\mathbb{X})$.

%

Let $h$ be a positive Borel measurable function
$h\dvtx \mathbb{R}^+\cup\{0\}\rightarrow\mathbb{R}^+$. Here
$h(\widetilde{\mu}(\mathbb{X}))$ is a tilting factor transforming
the total mass to a positive scalar. The law of the {\it tilted
completely random measure} (tilted {\sc crm}) is given by scaling
the law of the {\sc crm} by the tilting factor. To ensure that the
law of the tilted {\sc crm} is proper, the proportional constant of
the law, $E[h(\widetilde{\mu}(\mathbb{X}))]$, is required
to be finite, that is,
%
\begin{eqnarray}
\label{ConditionEhx} E\bigl[h\bigl(\widetilde{\mu}(\mathbb{X})\bigr)\bigr]<
\infty.
\end{eqnarray}
%
\begin{definition}\label{Deftiltedmeasure} Let $\widetilde{\mu}$ be a {\sc crm} defined on
$(\mathbb{X},\mathcal{B}(\mathbb{X}))$ in (\ref{CRM}). The {\sc crm}
$\widetilde{\mu}$ has a probability measure
$\mathcal{P}_{\widetilde{\mu}}$ on
$(\mathcal{M}_{\mathbb{X}},\mathcal{B}(\mathcal{M}_{\mathbb{X}}))$
and with the parameter measure $\nu$ that satisfies conditions
(\ref{condition1}) and (\ref{condition2}).
Let $h$ be a positive Borel measurable function on the non-negative
real line, that satisfies condition (\ref{ConditionEhx}). A tilted
{\sc
crm}, $\widetilde{\mu}_t$, defined on
$(\mathbb{X},\mathcal{B}(\mathbb{X}))$, has a probability measure
$\mathcal{P}_{\widetilde{\mu}_t}$ on
$(\mathcal{M}_{\mathbb{X}},\mathcal{B}(\mathcal{M}_{\mathbb{X}}))$
such that
%
\begin{eqnarray}
\label{tiltedlawsacledh}%
\mathcal{P}_{\widetilde{\mu}_t}(A):= 
\frac{1}{E[h(\widetilde{\mu}(\mathbb{X}))]}\int_{A}h\bigl(\mu(\mathbb {X})\bigr)
\mathcal{P}_{\widetilde{\mu}}(\mathrm{d}\mu), \qquad A\in\mathcal{B}(\mathcal {M}_{\mathbb{X}}),%
\end{eqnarray}
with the parameter measures $\nu$ and $h$.
\end{definition}
%
%
\begin{definition}\label{Deftiltedprobabilitymeasure} Let $\widetilde{\mu}_t$ be a tilted {\sc
crm} with the parameter measures $\nu$ and $h$ defined from
Definition \ref{Deftiltedmeasure}, a {\it normalized tilted} {\sc
crm} $\widetilde{G}_t$ is given by
%
\begin{eqnarray}
\label{tiltedrandomprobabilitymeasure}%
\widetilde{G}_t(A)= \widetilde{\mu}_t(A) / \widetilde{
\mu}_t(\mathbb {X}),\qquad %
A\in\mathcal{B}(\mathbb{X}),
\end{eqnarray}
on $(\mathbb{X},\mathcal{B}(\mathbb{X}))$. This normalized tilted
\textsc{crm} $\widetilde{G}_t$ is with the parameter measures $\nu$ and
$h$.
\end{definition}

Notice that when the function is a finite constant
($h(x)=\mathrm{constant}<\infty$), the normalized tilted \textsc{crm}
$\widetilde{G}_t$ is simply a normalized \textsc{crm},
which has been extensively studied by James \cite{James2002} and
James \textit{et al.} \cite{JLP2008}. Some special cases with various choices of function $h$
and \textsc{crm}s have been considered. James \cite{James2005b} considers
polynomially tilting the law of gamma process with
$h(\mu(\mathbb{X}))=\mu(\mathbb{X})^{-q}$.
Pitman \cite{Pitman2003,Pitman2006} constructs the Poisson--Dirichlet
process from polynomially tilting the law of positive stable
process. These are all interesting special cases covered by the class of
the normalized tilted \textsc{crm}s which will be further discussed in
Section \ref{SectionNRM-PET}.

The total mass is a key ingredient of both normalized \textsc{crm}s and
normalized tilted \textsc{crm}s. We consider the connection between these
two total masses and the general framework on a characterization of
the masses through the Laplace transform.
Let the total masses of \textsc{crm} and tilted \textsc{crm} be
$\widetilde{T}:=\widetilde{\mu}(\mathbb{X})$ and
$\widetilde{T}_t:=\widetilde{\mu}_t(\mathbb{X})$. Both
$\widetilde{\mu}$ and $\widetilde{\mu}_t$ are the mappings to the
positive line so that $\widetilde{T}>0$ and $\widetilde{T}_t>0$
and the laws of $\widetilde{T}$ and $\widetilde{T}_t$
are both absolutely continuous with respect to
Lebesgue measure. Their densities,
$f_{\widetilde{T}}(y)$ and $f_{\widetilde{T}_t}(y)$,
are related through the equality
%
\begin{eqnarray}
\label{densityT} f_{\widetilde{T}_t}(y) =\frac{h(y)f_{\widetilde{T}}(y)}{\int_{\mathbb{R}^+}h(y)
f_{\widetilde{T}}(y)\,\mathrm{d}y} .
\end{eqnarray}
The Laplace transform of the random variable $\widetilde{T}$ is given by
\begin{eqnarray*}
E \bigl[\mathrm{e}^{-\lambda\widetilde{T}} \bigr]=\int_{\mathbb{R}^{+}}\mathrm{e}^{-\lambda
y}f_{\widetilde{T}}(y)\,\mathrm{d}y%
=\mathrm{e}^{-\psi_{0}(\lambda)} ,
\end{eqnarray*}
where in general $\psi_{a}(b)$ is given by
%
\begin{eqnarray}
\label{psi} \psi_{a}(b)=\int_{\mathbb{R}^{+}\times\mathbb{X}}
\bigl(1-\mathrm{e}^{-b
s}\bigr)\mathrm{e}^{-as}\rho_x(\mathrm{d}s)\eta(\mathrm{d}x) ,\qquad a
\geq0 , b>0.%
\end{eqnarray}
In terms of the density of $\widetilde{T}$, the Laplace transform of
the random variable $\widetilde{T}_t$ could be seen as
%
\begin{eqnarray}
\label{LaplaceTt} E \bigl[\mathrm{e}^{-\lambda\widetilde{T}_t} \bigr]=%
\int
_{\mathbb{R}^{+}}\mathrm{e}^{-\lambda y}f_{\widetilde{T}_t}(y)\,\mathrm{d}y=%
\int
_{\mathbb{R}^{+}}\mathrm{e}^{-\lambda y}\frac{h(y)
f_{\widetilde{T}}(y)}{\int_{\mathbb{R}^+}h(s)
f_{\widetilde{T}}(s)\,\mathrm{d}s}\,\mathrm{d}y.
\end{eqnarray}
%
From (\ref{LaplaceTt}), it requires to specify $h$ to derive an
explicit form. In fact, it is clear that
the total masses are connected through the equality (\ref{densityT}) and this could be utilized to derive the distributional results of
the masses. A further extension on characterizing the random measures may
concern the Laplace transform of their functionals, specifically linear
functionals,
that play an important role in the studies of random measures.
One could see it as a
generalization from measuring a set (e.g., equation (\ref{CRM})) as an
indicator function to measuring a class of functions. Let
$BM^+(\mathbb{X})$ denote the class of positive Borel measurable
functions mapping from $\mathbb{X}$ to $\mathbb{R}^+\cup\{0\}$ and
all functions in this class vanish outside the bounded sets of
$\mathbb{X}$. Let $g$ be a function in $BM^+(\mathbb{X})$ and let
the functional defined as
$\widetilde{\mu}(g):=\int_{\mathbb{X}}g(x)\widetilde{\mu}(\mathrm{d}x)$. This
functional is also regarded as a Poisson functional since it has an
expression with respect to the Poisson process on $\mathbb{X}\times
\mathbb{R}^+$, that is
$\widetilde{\mu}(g)=\int_{\mathbb{R}^+\times\mathbb{X}}g(x)yN(\mathrm{d}x,\mathrm{d}y)$.
Then the properties of the \textsc{crm} $\widetilde{\mu}$ could be
derived from the Poisson process $N$.
For general discussion of functionals and
Laplace functionals, see Daley and Vere-Jones \cite{DV2008}, Chapters 9, 10. An early application to Bayesian non-parametric
statistics would be found in \cite{LW1989} (see also Dykstra and Laud \cite{DL1981}).

\section{The normalized random measure derived from the
polynomially and exponentially tilted law} \label{SectionNRM-PET}

Here we aim at showing the conjugacy property of the
normalized tilted \textsc{crm} $\widetilde{G}_t$ with a specific choice of $h$
(Definition \ref{Deftiltedprobabilitymeasure}). First, we take $h$ as follows
%
\begin{eqnarray}
\label{generalh} h(x)=h^{\prime}(x)x^{-q} ,\qquad q\geq0 ,
\end{eqnarray}
where $h^{\prime}$ is a positive measurable function
$h^{\prime}\dvtx \mathbb{R}^+\cup\{0\}\rightarrow\mathbb{R}^+$ and
satisfying condition (\ref{ConditionEhx}),
that is
$E[h^{\prime}(\widetilde\mu(\mathbb{X}))\widetilde\mu(\mathbb
{X})^{-q}]<\infty$,
but it is not depending on the scalar $q$.
Take the normalized tilted \textsc{crm} $\widetilde{G}_t$ be a prior
random probability measure in the Bayesian non-parametric content,
we can show that the posterior of $\widetilde{G}_t$ is belong to the
same class of random probability measure (Definition \ref{Deftiltedprobabilitymeasure}). So the conjugacy property of
$\widetilde{G}_t$ with the choice of $h$ (\ref{generalh}) is
immediately revealed.
We then consider the polynomially and
exponentially tilted law of the \textsc{crm} $\widetilde{\mu}_t$,
specifically, we take
%
\begin{eqnarray}
\label{PET} h^\prime(x)=\mathrm{e}^{-\gamma x},\qquad \gamma\geq0 ,
\end{eqnarray}
in (\ref{generalh}). This choice covers a rich class of random
probability measures and the interest in this choice is desirable.
The posterior analysis of this class of tilted \textsc{crm}s is given
after the conjugacy property has been shown.

In general, the law of the tilted \textsc{crm} $\widetilde{\mu}_t$ is
given by
\begin{eqnarray*}
\label{tiltedlaw} \mathcal{P}_{\widetilde{\mu}_t}(A)=%
\frac{1}%
{E[h^{\prime}(\widetilde{\mu}(
\mathbb{X}))\widetilde{\mu}(\mathbb{X})^{-q}]} \int
_{A} h^\prime\bigl(\mu(\mathbb{X})\bigr)\mu(
\mathbb{X})^{-q}%
\mathcal{P}_{\widetilde{\mu}}(\mathrm{d}\mu) ,\qquad A\in
\mathcal{B}(\mathcal {M}_\mathbb{X}) .%
\end{eqnarray*}
The parameters of both the tilted
\textsc{crm} $\widetilde{\mu}_t$ and normalized tilted \textsc{crm}
$\widetilde{G}_t$ are now $q\geq0$, $h^\prime$ and the intensity
measure $\nu$ only. In particular, if $h^\prime$ is chosen to be
$\mathrm{e}^{-\gamma x}$ in (\ref{generalh}), the parameters are then
$q\geq0$, $\gamma\geq0$ and the intensity measure $\nu$. The random
probability measures
in the class of normalized tilted \textsc{crm}s
are more general than those in the class of the normalized \textsc{crm}s (Regazzini \textit{et al.} \cite{RLP2003}; James \textit{et al.} \cite{JLP2008}); one could easily realize
that the tilted \textsc{crm} $\widetilde{\mu}_t$ (Definition \ref{Deftiltedmeasure}) is
not necessarily a \textsc{crm} and this could be seen as a generalization
of the class.

\begin{remark}\label{remark3exponentialtilting}
A referee pointed out that here the exponentially tilting, say
$h(x)=\mathrm{e}^{-\gamma x}$, is redundant as the exponentially tilting
operation on a \textsc{crm} leads to another \textsc{crm}. So taking
$h^\prime$ as in (\ref{PET}), such that $h(x)=\mathrm{e}^{-\gamma x}x^{-q}$,
operated with a \textsc{crm} is equivalent to taking $h^\prime$ as a
fixed finite constant, such that $h(x)=x^{-q}$, operated with
another \textsc{crm}.
\end{remark}
%
\begin{remark}\label{remark3crm}
Scaling operation on the law of Poisson random measures and \textsc{crm}s has been studied in James \cite{James2002}. James
\cite{James2002}, equation 70, has the same construction as in this work by taking
$h^\prime$ to be a fixed finite constant in (\ref{generalh}), or to
be (\ref{PET}) from Definition \ref{Deftiltedmeasure}.
\end{remark}

\subsection{The posterior law and structural conjugacy}
Consider a sequence of exchangeable random elements taking values in
$\mathbb{X}$. These random variables are assumed to be conditionally
independent and identically distributed given the normalized tilted
\textsc{crm} $\widetilde{G}_t$ (Definition \ref{Deftiltedprobabilitymeasure}) such that for every integer $n\geq1$
%
\begin{eqnarray}
\label{likelihood} \mathbb{P}(X_1\in B_1,
\ldots,X_n\in B_n|\widetilde{G}_t)=%
\prod_{i=1}^n \widetilde{G}_t(B_i)=
\prod_{i=1}^n \frac{\widetilde{\mu}_t(B_i)}{\widetilde{\mu}_t(\mathbb
{X})}
,\qquad%
B_1,\ldots,B_n\in\mathcal{B}(\mathbb{X}) ,\quad
\end{eqnarray}
where $\widetilde{G}_t$, or equivalently $\widetilde{\mu}_t$, is
regarded as a parameter. Then, (\ref{likelihood}) is the
``likelihood'' with the parameter $\widetilde{\mu}_t$. Let
$\mathcal{P}_{\widetilde{\mu}_t}^{(\mathbf{X}_n)}$ indicate the
posterior distribution of $\widetilde{\mu}_t$, namely the
distribution of $\widetilde{\mu}_t$ conditional on $\mathbf{X}_n$,
that is
\begin{eqnarray*}
\mathcal{P}_{\widetilde{\mu}_t}^{(\mathbf{X}_n)}(A):= \frac{\mathcal{P}_{\widetilde{\mu}_t,\mathbf{X}_n}(A,B_1,\ldots,B_n)} {
\mathcal{P}_{\widetilde{\mu}_t,\mathbf{X}_n}(\mathcal{M}_\mathbb
{X},B_1,\ldots,B_n)} ,\qquad
A\in\mathcal{B}(\mathcal{M}_\mathbb{X}) ,%
B_1,\ldots,B_n\in\mathcal{B}(\mathbb{X}) ,
\end{eqnarray*}
where $\mathcal{P}_{\widetilde{\mu}_t,\mathbf{X}_n}$ represents the
joint distribution of $(\widetilde{\mu}_t,\mathbf{X}_n)$ such that
%
\begin{eqnarray}
\label{jointlaw} 
\mathcal{P}_{\widetilde{\mu}_t,\mathbf{X}_n}(A,B_1,
\ldots,B_n)=%
\int
_{A}\prod_{i=1}^n
\frac{\mu(B_i)}{\mu(\mathbb{X})}\mathcal {P}_{\widetilde{\mu}_t}(\mathrm{d}\mu) , \qquad
A\in
\mathcal{B}(\mathcal{M}_\mathbb{X}) , B_1,
\ldots,B_n\in\mathcal{B}(\mathbb{X}) . \qquad
\end{eqnarray}
Taking into consideration of the above assumption of
conditional independence (\ref{likelihood}), the joint distribution
$\mathcal{P}_{\widetilde{\mu}_t,\mathbf{X}_n}$of $(\widetilde{\mu
}_t,\mathbf{X}_n)$ is defined on the space
$(\mathcal{M}_\mathbb{X}\times\mathbb{X}^n,\mathcal{B}(\mathcal
{M}_\mathbb{X})\otimes\mathcal{B}(\mathbb{X})^{n})$
for $n=1,2,\ldots$\,. Now, considering Definition \ref{Deftiltedmeasure} of
$\mathcal{P}_{\widetilde{\mu}_t}$, one obtains
%
\begin{eqnarray}
\label{expression1} &&\int_{A} \prod
_{i=1}^n \frac{\mu(B_i)}{\mu(\mathbb{X})}\mathcal
{P}_{\widetilde{\mu}_t}(\mathrm{d}\mu)=%
\frac{1}%
{E
[h^\prime(\widetilde{\mu}(\mathbb{X}))\widetilde{\mu}(
\mathbb{X})^{-q}]} \int_{A}
\frac{h^\prime(\mu(\mathbb{X}))}{\mu(\mathbb{X})^{n+q}}%
\Biggl(\prod_{i=1}^n
\mu(B_i) \Biggr) \mathcal{P}_{\widetilde{\mu}}(\mathrm{d}\mu) ,
\quad\\
&& \quad
A\in\mathcal{B}(\mathcal{M}_\mathbb{X}) ,
B_1,\ldots,B_n\in\mathcal{B}(\mathbb{X}) .
\nonumber
\end{eqnarray}
This shows the key element needed for deriving the posterior law of
$\widetilde{\mu}_t$
is the law of the \textsc{crm} $\mathcal{P}_{\widetilde{\mu}}$ as seen in
the right-hand side
of (\ref{expression1}). The usual technique to derive the posterior of
$\widetilde{\mu}_t$ requires application of change of
measure or disintegration. So, the major task is to
apply change of measure updating the law $\mathcal{P}_{\widetilde{\mu
}}$ with the information,
namely $(\prod_{i=1}^n \mu(B_i))$,
$h^\prime({\mu}(\mathbb{X}))$ and $\mu(\mathbb{X})^{-(n+q)}$,
to a posterior law, $\mathcal{P}_{\widetilde{\mu}}^{(\mathbf{X}_n)}$.
Dealing with the term $(\prod_{i=1}^n \mu(B_i))$ could be simply adopted
by the standard arguments (James \cite{James2002,James2005a,James2005b}).
The next term $h^\prime({\mu}(\mathbb{X}))$ should be chosen explicitly
to proceed.
When $h^\prime({\mu}(\mathbb{X}))=\mathrm{e}^{-\gamma{\mu}(\mathbb{X})}$, change
of measure
only involves the Laplace transform and doesn't cost much effort. Eventually,
dealing with the term
$\mu(\mathbb{X})^{-(n+q)}$ could be somehow challenging. James \cite{James2002} (see also
James \textit{et al.} \cite{JLP2008}) introduces an augmentation approach that allows us to
proceed further in particular for this polynomial term $\mu(\mathbb
{X})^{-(n+q)}$.
This leads to the analysis on the posteriors of the normalized tilted
\textsc{crm}s.

James \cite{James2002}'s approach makes use of the gamma identity and
introduces an augmented
variable. We now address the role of the augmented variable. Here
the well-known gamma identity is given by
%
\begin{eqnarray}
\label{gammaidentity} \frac{1}{b^a}=\frac{1}{\Gamma(a)} \int
_{\mathbb{R}^{+}}%
\mathrm{e}^{-bu}u^{a-1}\,\mathrm{d}u ,\qquad a,b\in
\mathbb{R}^+.%
\end{eqnarray}
Then take the term $\mu(\mathbb{X})^{-(n+q)}$ in (\ref{expression1}) as
$b^{-a}$ in (\ref{gammaidentity}),
term involving $\mu$ becomes tractable, positioned as an exponent. The
integral (numerator) appearing in the right-hand side
(\ref{expression1}) can be rewritten as
%
\begin{eqnarray}
\label{expression2} 
&& 
\int_{A\times\mathbb{R}^{+}}%
h^\prime(y)%
\biggl[\int_{\mathbb{R}^{+}}
\frac{u^{n+q-1}\mathrm{e}^{-u y}}{\Gamma(n+q)}\,\mathrm{d}u \biggr] 
\Biggl(\prod
_{i=1}^n \mu(B_i) \Biggr) \mathbb{P}(
\widetilde{\mu}\in \mathrm{d}\mu|\widetilde{T}=y)f_{\widetilde{T}}(y)\,\mathrm{d}y , 
\\
&& \quad
A\in\mathcal{B}(\mathcal{M}_\mathbb{X}) ,
B_1,\ldots,B_n\in\mathcal{B}(\mathbb{X}) .
\nonumber
\end{eqnarray}
where $\mathbb{P}(\widetilde{\mu}\in d\mu|\widetilde{T}=y)$ denotes
the conditional distribution of $\widetilde{\mu}$ given its' total
mass $\widetilde{T}=\widetilde{\mu}(\mathbb{X})$.
Rewriting the last expression by replacing $f_{\widetilde{T}}$ with
its expression in term $f_{\widetilde{T}_t}$ (\ref{densityT}), from
(\ref{expression1}) and (\ref{expression2}) one finally gets
%
\begin{eqnarray}
&&\int_{A} \prod_{i=1}^n
\frac{\mu(B_i)}{\mu(\mathbb{X})}\mathcal {P}_{\widetilde{\mu}_t}(\mathrm{d}\mu)%
=%
\int_{A\times\mathbb{R}^{+}\times\mathbb{R}^{+}}%
\Biggl(\prod
_{i=1}^n \frac{\mu(B_i)}{y} \Biggr) \mathbb{P}(
\widetilde{\mu}\in \mathrm{d}\mu|\widetilde{T}=y)%
\ell_n(y,u)\,\mathrm{d}y\,\mathrm{d}u ,
\quad\\
&&\quad
A\in\mathcal{B}(\mathcal{M}_\mathbb{X}) ,
B_1,\ldots,B_n\in\mathcal{B}(\mathbb{X}) .
\nonumber
\end{eqnarray}
where
%
\begin{eqnarray}
\ell_n(y,u)=\frac{u^{n+q-1}\mathrm{e}^{-u
y}y^{n+q}}{\Gamma(n+q)}f_{\widetilde{T}_t}(y) ,\qquad
u>0
, y>0 .%
\end{eqnarray}
Here $\ell_n$ is a joint probability density on
$\mathbb{R}^+\times\mathbb{R}^+$. Without loss of generality, we
assume that $(\Omega,\mathfrak{F},\textsf{P})$ is large enough to
support a sequence $U_n$ of random numbers such that the distribution of
$(U_n,\widetilde{T}_t)$ admits $\ell_n$ as density function. Then,
it follows that
%
\begin{eqnarray}
\label{Unmarginaldistribution} f_{U_n}(u)=%
\int
_{\mathbb{R}^{+}}\ell_n(y,u)\,\mathrm{d}y= \frac{u^{n+q-1}}{\Gamma(n+q)}%
\int_{\mathbb{R}^{+}}y^{n+q}\mathrm{e}^{-uy}%
f_{\widetilde{T}_t}(y)\,\mathrm{d}y ,\qquad u>0 .%
\end{eqnarray}
Here $f_{U_n}$ is a probability density for $U_n$.
In view of these elementary developments, one can disintegrate the
law of $(\widetilde{\mu}_t,\mathbf{X}_n)$ as follows
%
\begin{eqnarray}
\label{expression3} && \int_{A} \prod
_{i=1}^n \frac{\mu(B_i)}{\mu(\mathbb{X})}\mathcal
{P}_{\widetilde{\mu}_t}(\mathrm{d}\mu)%
=\int
_{A\times\mathbb{R}^{+}\times\mathbb{R}^{+}}%
\Biggl(\prod_{i=1}^n
\frac{\mu(B_i)}{y} \Biggr) 
\nonumber\\
&&\hspace*{141pt}{}\times\mathbb{P}(\widetilde{\mu}\in \mathrm{d}\mu|
\widetilde{T}=y)%
f_{\widetilde{T}|U_n}(y|u)%
f_{U_n}(u)%
\,\mathrm{d}u\,\mathrm{d}y ,\quad
\\
&&\quad
A\in\mathcal{B}(\mathcal{M}_\mathbb{X}) ,
B_1,\ldots,B_n\in\mathcal{B}(\mathbb{X}) ,
\nonumber
\end{eqnarray}
where
\begin{eqnarray*}
f_{\widetilde{T}|U_n}(y|u)=\frac{\ell_n(u,y)}{\int_{\mathbb{R}^{+}}\ell_n(u,y)\,\mathrm{d}y}=
\frac{y^{n+q}\mathrm{e}^{-uy}}{E[\widetilde{T}_t^{n+q}\mathrm{e}^{-u\widetilde{T}_t}]}%
f_{\widetilde{T}_t}(y)%
,\qquad y>0 .
\end{eqnarray*}
This disintegration (\ref{expression3}) shows the role of $U_n$ as
an augmented variable.
%
\begin{remark}
The representation of the density (\ref{Unmarginaldistribution})
suggests that $U_n$ has the gamma distribution with a random scale
$\widetilde{T}_t$ which has the density $f_{\widetilde{T}_t}$. Given
$\widetilde{T}_t$, $U_n$ has the gamma distribution with parameters
$(n+q,\widetilde{T}_t)$. The product of the random variables,
$U_n\widetilde{T}_t$, has the gamma distribution with parameters
$(n+q,1)$ and independent of $\widetilde{T}_t$.
\end{remark}

Before proceeding to the posterior distribution, additional notations
are introduced. The normalized tilted \textsc{crm} $\widetilde{G}_t$ is
almost surely discrete. A random sample $\mathbf{X}_n$ of $\widetilde{G}_t$
usually contains ties. We can always express $\mathbf{X}_n$ by two elements,
namely a partition and distinct values.
Here $\mathbf{p}_n$ is a partition of the
integers $ \{1,\ldots,n \}$ that are the indices of $\mathbf
{X}_n$ and
$\mathbf{Y}_{\aleph(\mathbf{p}_n)}= \{Y_1,\ldots,Y_{\aleph(\mathbf
{p}_n)} \}$
represents the distinct values of $\mathbf{X}_n$.
The partition
$\mathbf{p}_n$ locates the distinct values from $\mathbf{X}_n$ to
$\mathbf{Y}_{\aleph(\mathbf{p}_n)}$ or vice versa. As a result, we have
the following equivalent representations
%
\begin{eqnarray}
\label{decomp} \mathbf{X}_n= \{X_1,\ldots,X_n
\}= \{Y_1,\ldots,Y_{\aleph
(\mathbf{p}_n)},\mathbf{p}_n \}= \{
\mathbf{Y}_{\aleph(\mathbf
{p}_n)},\mathbf{p}_n \} .
\end{eqnarray}
A partition $\mathbf{p}_n$ contains
$\aleph(\mathbf{p}_n)$ cells (known as clusters), that is
$\mathbf{p}_n= \{C_{1},\ldots,C_{\aleph(\mathbf{p}_n)} \}$.
Each cell $C_{k}$ contains the indices of a subset of
$\mathbf{X}_n$,
namely the unique values $Y_k$
such that $C_{k}=  \{i\dvtx X_i=Y_k,i=1,\ldots,n \}$ for
$k=1,\ldots,\aleph(\mathbf{p}_n)$. The number of elements in the
cell $k$, $C_k$, of the partition is indicated by $n_k$, for
$k=1,\ldots,\aleph(\mathbf{p}_n)$, so that
$\sum_{k=1}^{\aleph(\mathbf{p}_n)} n_{k}=n$. Therefore, the union of
all cells is the set of all $n$ integers,
$\bigcup_{k=1}^{\aleph(\mathbf{p}_n)} C_{k} =
\{1,\ldots,n \}$ and all cells are pairwise mutually
exclusive, that is $C_{k} \cap C_{k^{\prime}} = \emptyset$ where
$k\neq k^\prime$ for $k,k^{\prime} =1,\ldots,\aleph(\mathbf{p}_n)$.
This partition representation is commonly used in Bayesian
non-parametric literature
(see Lo \cite{Lo1984}; Lo and Weng \cite{LW1989}; James \cite{James2002,James2005a,James2005b}) since it
well describes
the variates generated from those random probability measures and
is also useful in expressing the marginal distribution of $\mathbf{X}_n$.

Lijoi and Pr\"unster \cite{LP2010} describe the concept of structural conjugacy. A
random probability measure, say $G$, is a structurally conjugate
random probability measure if the resulting posterior law of $G$
given $\mathbf{X}_n= \{X_1,\ldots,X_n \}$, has the same
structure. {\it Neutral to the right} process (Doksum \cite{Doksum1974})
is one of the classes that has this property. In the present work,
we show that the normalized tilted \textsc{crm}s in Definition
\ref{Deftiltedmeasure} are also structurally conjugate. Here we
follow James \textit{et al.} \cite{JLP2008},
write $\widetilde{G}_t^{(\mathbf{X}_n)}$ as the posterior normalized
tilted \textsc{crm}.
Theorem \ref{thm1posteriorG} shows that the
normalized tilted \textsc{crm}s in Definition \ref{Deftiltedprobabilitymeasure} have the conjugate property, that is, both
$\widetilde{G}_t$ and $\widetilde{G}_t^{(\mathbf{X}_n)}$ are in the
same class.
%
\begin{theorem}\label{thm1posteriorG}
Let
$\widetilde{G}_t=\widetilde\mu_t /\widetilde\mu_t(\mathbb{X})$ be
a normalized tilted \textsc{crm} (Definition \ref{Deftiltedprobabilitymeasure}) defined on $\mathbb{X}$ with $h$ as in
(\ref{generalh}), that is $h(x)=h^\prime(x)x^{-q}$. The parameters
of this normalized tilted \textsc{crm} are given by the Borel
measurable function $h^\prime$, the measure $\nu$, and the scalar
$q\geq0$. With the prior measure $\widetilde{G}_t$ and the likelihood
(\ref{likelihood}),
the posterior measure, namely $\widetilde{G}_t^{(\mathbf{X}_n)}$, has
the same distribution of a normalized tilted
\textsc{crm} such that
%
\begin{eqnarray}
\widetilde{G}_t^{(\mathbf{X}_n)}
\overset{d}
{=}%
\widetilde\mu_t^{(\mathbf{X}_n)}
/%
\widetilde\mu_t^{(\mathbf{X}_n)}(\mathbb{X}) ,%
\end{eqnarray}
where
\begin{enumerate}
\item$\widetilde\mu_t^{(\mathbf{X}_n)}$ has the law %
\begin{eqnarray}
\label{tiltedposterior} \mathcal{P}_{\widetilde\mu_t}^{(\mathbf{X}_n)}(A) =
\frac{1}%
{E[h^\prime(\widetilde
\mu^{(\mathbf{X}_n)}(\mathbb{X})) \widetilde\mu^{(\mathbf{X}_n)}(
\mathbb{X})^{-(n+q)}]}\int_A h^\prime
\bigl(\mu(\mathbb{X})\bigr) \mu(\mathbb{X})^{-(n+q)} \mathcal{P}_{\widetilde\mu}^{(\mathbf{X}_n)}(\mathrm{d}\mu)
\nonumber
\end{eqnarray}
for $A\in\mathcal{B}(\mathcal{M}_\mathbb{X})$.
\item$\widetilde\mu^{(\mathbf{X}_n)}$ is a normalized tilted
\textsc{crm} whose law $\mathcal{P}_{\widetilde\mu}^{(\mathbf{X}_n)}$,
has the same law as
%
\begin{eqnarray}
\widetilde\mu+\sum_{k=1}^{\aleph(\mathbf{p}_n)}J_k^{(\mathbf
{X}_n)}
\delta_{Y_k},
\end{eqnarray}
where $\delta_{a}$ is the Dirac delta function evaluated at $a$,
$\widetilde\mu$ is a \textsc{crm} with law
$\mathcal{P}_{\widetilde\mu}$ and intensity measure $\nu$,
$ \{Y_1,\ldots,Y_{\aleph(\mathbf{p}_n)} \}$ is a sequence of
fixed points of discontinuity, and
$ \{J_{1}^{(\mathbf{X}_n)},\ldots,J_{\aleph(\mathbf{p}_n)}^{(\mathbf
{X}_n)} \}$
are the corresponding jumps.
\item Conditional on $\mathbf{X}_n$, each jump $J_{k}^{(\mathbf{X}_n)}$
has the
conditional distribution
\begin{eqnarray*}
\mathbb{P}\bigl(J_k^{(\mathbf{X}_n)}\in \mathrm{d}s\bigr)=%
\frac{s^{n_k}\rho_{Y_k}(\mathrm{d}s)}%
{\int_{\mathbb{R}^{+}}s^{n_k}
\rho_{Y_k}(\mathrm{d}s)} , \qquad%
k=1,\ldots,\aleph(\mathbf{p}_n)
.
\end{eqnarray*}
\item Conditional on $\mathbf{X}_n$, $\widetilde{\mu}$ and
$ \{J_{1}^{(\mathbf{X}_n)},\ldots,J_{\aleph(\mathbf{p}_n)}^{(\mathbf
{X}_n)} \}$
are independent.
\end{enumerate}
\end{theorem}
\begin{pf}
See Appendix \ref{thm1proof}.
\end{pf}

In Theorem \ref{thm1posteriorG}, the posterior law of
$\widetilde{\mu}_t$,
namely $\widetilde{\mu}_t^{(\mathbf{X}_n)}$, is shown to have the tilted
law of the \textsc{crm} $\widetilde\mu^{(\mathbf{X}_n)}$, which has the law
$\mathcal{P}_{\widetilde\mu}^{(\mathbf{X}_n)}$, where the tilting
factor is updated to
$h^\prime(\widetilde{\mu}^{(\mathbf{X}_n)}(\mathbb{X}))
\widetilde{\mu}^{(\mathbf{X}_n)}(\mathbb{X})^{-(n+q)}$. Under the
normalization of the process, the posterior of the random
probability measure becomes
$\widetilde{G}_t^{(\mathbf{X}_n)}\overset{d}=\widetilde{\mu}_t^{(\mathbf
{X}_n)} /\widetilde\mu_t^{(\mathbf{X}_n)}(\mathbb{X})$.
This confirms the conjugate property of the normalized tilted \textsc{crm}.
%
\begin{remark}\label{remark3thm1}
Following Remark \ref{remark3crm}, the posterior distribution
discussed in James \cite{James2002}, Chapter~5, over the scaling operation
on \textsc{crm}s has been established in James \cite{James2002}, Corollary 5.1. This result is a version of Theorem \ref{thm1posteriorG}.
Furthermore, James \cite{James2002}, Theorem 5.1, also
supplies the posterior law of the corresponding Poisson random
measure. In particular, putting statements i and ii of
James \cite{James2002}, Corollary 5.1, and statements i and ii of
James \cite{James2002}, Proposition 5.2, together is equivalent to Theorem
\ref{thm1posteriorG}.
\end{remark}

A special case of the tilted \textsc{crm} (Definition \ref{Deftiltedmeasure}) that takes $h(x)=\mathrm{e}^{-\gamma x}x^{-q}$ is of interest. The
following theorem, Theorem \ref{thm2conditionalposteriorG},
describes an augmented posterior law of $\widetilde{G}_t$ and
$\widetilde{\mu}_t$, denoted by
$\widetilde{G}_t^{(U_n,\mathbf{X}_n)}$ and
$\widetilde{\mu}_t^{(U_n,\mathbf{X}_n)}$, that are the conditional
laws of $\widetilde{G}_t$ and $\widetilde{\mu}_t$, respectively, given
both $U_n$ and $\mathbf{X}_n$.
%
\begin{theorem}\label{thm2conditionalposteriorG}
Let
$\widetilde{G}_t=\widetilde\mu_t /\widetilde\mu_t(\mathbb{X})$ be
a normalized tilted \textsc{crm} (Definition \ref{Deftiltedprobabilitymeasure}) on $\mathbb{X}$ with $h$ as in (\ref{generalh}) and $h^\prime$ as in (\ref{PET}), that is $h(x)=\mathrm{e}^{-\gamma
x}x^{-q}$. The parameters of this normalized tilted \textsc{crm} are
given by the measure $\nu$, and the scalers $q\geq0$ and
$\gamma\geq0$. With the prior measure $\widetilde{G}_t$ and the
likelihood (\ref{likelihood}),
the conditional posterior measure, namely
$\widetilde{G}_t^{(U_n,\mathbf{X}_n)}$, given $U_n$ and $\mathbf
{X}_n= \{X_1,\ldots,X_n \}$
has the same distribution of a normalized
\textsc{crm} such that
\begin{eqnarray*}
\widetilde{G}_t^{(U_n,\mathbf{X}_n)}
\overset{d} {=}
\widetilde{\mu}^{(U_n,\mathbf{X}_n)}
/\widetilde{\mu}^{(U_n,\mathbf{X}_n)}(
\mathbb{X}) ,%
\end{eqnarray*}
where
\begin{enumerate}
\item$\widetilde{\mu}^{(U_n,\mathbf{X}_n)}$ is a \textsc{crm} with the
same law
as
%
\begin{eqnarray}
\label{mugivenU} \widetilde{\mu}^{(U_n)} +\sum
_{k=1}^{\aleph(\mathbf{p}_{n})} J_{k}^{(U_n,\mathbf{X}_n)}
\delta_{Y_k} ,
\end{eqnarray}
where $\widetilde{\mu}^{(U_n)}$ is a \textsc{crm} with intensity
$\nu(\mathrm{d}x,\mathrm{d}s)=\eta(\mathrm{d}x)\times \mathrm{e}^{-s(\gamma+U_n)}\rho_x(\mathrm{d}s)$,
$ \{Y_1,\ldots,\break Y_{\aleph(\mathbf{p}_n)} \}$ is a sequence of
fixed points of discontinuity, and
$ \{J_{1}^{(U_n,\mathbf{X}_n)},\ldots,J_{\aleph(\mathbf
{p}_n)}^{(U_n,\mathbf{X}_n)} \}$
are the corresponding jumps. %
\item Conditional on $U_n$ and $\mathbf{X}_n$, each jump
$J_{k}^{(U_n,\mathbf{X}_n)}$ has the
conditional density
\begin{eqnarray*}
\mathbb{P}\bigl(J_k^{(U_n,\mathbf{X}_n)}\in \mathrm{d}s\bigr)= \frac{s^{n_k}\mathrm{e}^{-s(\gamma+U_n)}\rho_{Y_k}(\mathrm{d}s)}{\int_{\mathbb{R}^{+}}
s^{n_k}\mathrm{e}^{-s(\gamma+U_n)}\rho_{Y_k}(\mathrm{d}s)}
,\qquad k=1,\ldots,\aleph(\mathbf{p}_n) .
\end{eqnarray*}
\item Conditional on $U_n$ and $\mathbf{X}_n$, $\mu^{(U_n)}$ and
$ \{J_{1}^{(U_n,\mathbf{X}_n)},\ldots,J_{\aleph(\mathbf
{p}_n)}^{(U_n,\mathbf{X}_n)} \}$
are independent.
\item Conditional on $\mathbf{X}_n$, $U_n$ has the density
%
\begin{eqnarray}
\label{UgivenX} f_{U_{n}\vert\mathbf{X}_{n}}(u)=%
\frac{\mathrm{e}^{-\psi_{0}(\gamma+u)}%
\prod_{k=1}^{\aleph(\mathbf{p}_{n})}%
\tau_{n_{k},Y_{k}}(\gamma+u)u^{n+q-1}}%
{
\int_{\mathbb{R}^{+}}%
\mathrm{e}^{-\psi_{0}(\gamma+u)}%
\prod
_{k^{\prime}=1}^{\aleph(\mathbf{p}_{n})}%
\tau_{n_{k^{\prime}},Y_{k^{\prime}}}(
\gamma+u)u^{n+q-1}\,\mathrm{d}u} ,%
\end{eqnarray}
where
%
\begin{eqnarray}
\label{tau} \tau_{m,z}(a)=\int_{\mathbb{R}^{+}}s^{m}\mathrm{e}^{-sa}
\rho_z(\mathrm{d}s) ,\qquad m>0 , z\in \mathbb{X} , a\geq0 ,%
\end{eqnarray}
and $\psi_a(b)$ for $a\geq0$ and $b\geq0$ is defined in (\ref{psi}).
\end{enumerate}
\end{theorem}
\begin{pf} See Appendix \ref{thm2proof}.
\end{pf}

Theorem \ref{thm2conditionalposteriorG} shows that the
augmented posterior random probability measure
$\widetilde{G}_t^{(U_n,\mathbf{X}_n)}$ has the same distribution of
a normalized \textsc{crm} and further also provides the posterior
distribution of the augmented variables $U_n$,
$f_{U_{n}\vert\mathbf{X}_{n}}$. Combining these two yields the
posterior random probability measure
$\widetilde{G}_t^{(\mathbf{X}_n)}$. This is certainly useful in
applications of Bayesian non-parametric. For example, this theorem
could be useful in simulating the posterior normalized tilted \textsc{crm}
$\widetilde{G}_t^{(\mathbf{X}_n)}$ that is desirable in some
applications.\looseness=1
%
\begin{remark}\label{remark3thm2}
The posterior distribution of the tilted \textsc{crm}, namely
$\widetilde{\mu}^{(\mathbf{X}_n)}_t$, can be achieved by mixing
the law in (\ref{mugivenU}) (Theorem \ref{thm2conditionalposteriorG}, statement 1) over the distribution of $U_n$ given $\mathbf{X}_n$
in (\ref{UgivenX}) (Theorem \ref{thm2conditionalposteriorG},
statement 4). This is equivalent to the law of the tilted \textsc{crm}
in statement i of James \cite{James2002}, Corollary 5.1.
\end{remark}
%

\subsection{Generalized Blackwell and MacQueen P\'{o}lya urn sampling
scheme and marginal distribution of partitions}
\label{SectionUrnandmarginal} Blackwell and MacQueen \cite{BM1973} first introduce the
P\'{o}lya urn sampling scheme for the Dirichlet process and this
scheme can be utilized to generate random sequences from the
Dirichlet process. This employs so called the {\it
Blackwell--MacQueen P\'olya} (\textsc{bmp}) urn formula, that is the
predictive distribution of the Dirichlet process random sequences.
Lo \cite{Lo1991} also shows that this \textsc{bmp} urn formula can be
used to characterize the Dirichlet process. James \textit{et al.} \cite{JLP2008}
generalize the Blackwell--MacQueen P\'olya (\textsc{bmp}) urn formula
for the normalized \textsc{crm}s, namely {\it Generalized
Blackwell--MacQueen P\'olya} (\textsc{gbmp}) urn formula. A further
generalization will be considered for the normalized tilted \textsc{crm} in this section.

We consider a normalized tilted \textsc{crm} $\widetilde G_t$ with
$h(x)=\mathrm{e}^{-\gamma x}x^{-q}$ in Definition \ref{Deftiltedprobabilitymeasure}. Here the \textsc{gbmp} urn formula for this
normalized tilted \textsc{crm} will be presented under two
formulations, namely the unconditional and the conditional \textsc{gbmp} urn formulas for the normalized tilted \textsc{crm}. The
unconditional \textsc{gbmp} urn formula is simply the predictive
distribution, $\mathbb{P} \{X_{n+1}\in \mathrm{d}x \vert
\mathbf{X}_{n} \}$ where $\{X_{n+1},\mathbf{X}_n\}$ is the random
sequence drawn from the normalized tilted \textsc{crm} $\widetilde
G_t$. The conditional \textsc{gbmp} urn formula is the augmented
version of \textsc{gbmp} urn formula derived from Theorem \ref{thm2conditionalposteriorG}. An impression directly comes to the mind
is that the conditional urn formula, namely
$\mathbb{P} \{X_{n+1}\in \mathrm{d}x\vert\mathbf{X}_{n},U_n \}$ could
be derived according to the predictive distribution
$\mathbb{P} \{X_{n+1}\in \mathrm{d}x\vert\mathbf{X}_{n} \}=
E[\mathbb{P} \{X_{n+1}\in \mathrm{d}x\vert\mathbf{X}_{n},U_n \}\vert
\mathbf{X}_{n}]$ with respect to the distribution $U_n$ given
$\mathbf{X}_{n}$, $f_{U_n\vert
\mathbf{X}_{n}}$ (\ref{UgivenX}). 
However, the term $\mathbb{P} \{X_{n+1}\in \mathrm{d}x\vert
\mathbf{X}_{n},U_n \}$ is not necessarily a proper distribution.
A rescaling over both $\mathbb{P} \{X_{n+1}\in \mathrm{d}x\vert
\mathbf{X}_{n},U_n \}$ and $f_{U_n\vert
\mathbf{X}_{n}}$ with a factor involving terms like $U_n$ seems to be
needed. This leads to a new variable $\widetilde{U}_n$ introduced and
this new variable has a tilted density of $U_n$.
So, we obtain a proper distribution $\mathbb{P} \{X_{n+1}\in
\mathrm{d}x\vert\mathbf{X}_{n},\widetilde{U}_n \}$ and also the conditional
\textsc{gbmp} urn formula is established by mixing over
$\widetilde{U}_n$.
%
\begin{proposition}\label{prop1BMUrnpredictivedistribution}
Let
$\widetilde{G}_t=\widetilde\mu_t /\widetilde\mu_t(\mathbb{X})$ be
a normalized tilted \textsc{crm} (Definition \ref{Deftiltedprobabilitymeasure}) defined on $\mathbb{X}$ with $h$ as in
(\ref{generalh}) and $h^\prime$ as in (\ref{PET}), that is
$h(x)=\mathrm{e}^{-\gamma x}x^{-q}$. The parameters of this normalized tilted
\textsc{crm} are given by the measure $\nu$, and the scalers $q\geq0$
and $\gamma\geq0$. Then,
\begin{enumerate}
\item The predictive distribution for $X_{n+1}$ given
$\mathbf{X}_n= \{X_1,\ldots,X_n \}$ is given by
%
\begin{eqnarray}
\label{tiltedCRMBM} \mathbb{P} \{X_{n+1}\in \mathrm{d}x\vert
\mathbf{X}_{n} \} =%
\frac{\omega_{n,\aleph(\mathbf{p}_n)+1}(x)}{\phi(\mathbf{X}_n)}\eta(\mathrm{d}x)%
+
\sum_{k=1}^{\aleph(\mathbf{p}_n)}%
\frac{\omega_{n,k}(Y_k)}{\phi(\mathbf{X}_n)}%
\delta_{Y_k}(\mathrm{d}x) ,
\end{eqnarray}
where
\begin{eqnarray*}
\phi(\mathbf{X}_n)&=&\int
_{\mathbb{X}}\omega_{n,\aleph(\mathbf
{p}_n)+1}(x)\eta(\mathrm{d}x)+\sum
_{k=1}^{\aleph(\mathbf{p}_n)}\omega_{n,k}(Y_k) ,
\\
\omega_{n,\aleph(\mathbf{p}_n)+1}(x)&=&\int_{\mathbb{R}^{+}}u
\tau_{1,x}(\gamma+u)f_{U_n\vert\mathbf{X}_n}(u)\,\mathrm{d}u \quad \mbox{and}
\\
\omega_{n,k}(Y_k)&=&\int_{\mathbb{R}^{+}}u
\frac{\tau_{n_k+1,Y_k}(\gamma
+u)}{\tau_{n_k,Y_k}(\gamma+u)}f_{U_n\vert\mathbf{X}_n}(u)\,\mathrm{d}u , \qquad%
k=1,\ldots,\aleph(
\mathbf{p}_n) .
\end{eqnarray*}
\item Conditional on $\widetilde{U}_n$, the predictive distribution for
$X_{n+1}$ given
$\mathbf{X}_n= \{X_1,\ldots,X_n \}$ is given by
%
\begin{eqnarray}
\label{ConditionalBMUrn} \mathbb{P} \{X_{n+1}\in \mathrm{d}x\vert
\mathbf{X}_{n},\widetilde{U}_n \} =%
\frac{\omega_{n,\aleph(\mathbf{p}_n)+1}(\widetilde{U}_n, x)}{\phi
(\widetilde{U}_n, \mathbf{X}_n)}%
\eta(\mathrm{d}x)+%
\sum
_{k=1}^{\aleph(\mathbf{p}_n)}%
\frac{\omega_{n,k}(\widetilde{U}_n, Y_k)}{\phi(\widetilde{U}_n, \mathbf
{X}_n)}%
\delta_{Y_k}(\mathrm{d}x) ,
\qquad\end{eqnarray}
where
\begin{eqnarray*}
\phi(u,\mathbf{X}_n) &=& \int_{\mathbb{X}}
\omega_{n,\aleph(\mathbf{p}_n)+1}(u,x)\eta(\mathrm{d}x)+\sum_{k=1}^{\aleph(\mathbf{p}_n)}
\omega_{n,k}(u,Y_k) ,
\nonumber
\\
\omega_{n,\aleph(\mathbf{p}_n)+1}(u,x)&=&u\tau_{1,x}(\gamma+u) ,
\nonumber
\\
\omega_{n,k}(u,Y_k)&=&u \frac{\tau_{n_k+1,Y_k}(\gamma+u)}{\tau_{n_k,Y_k}(\gamma+u)}
,\qquad%
k=1,\ldots,\aleph(\mathbf{p}_n) \quad \mbox{and}
\nonumber
\\
\label{UprimegivenX} f_{\widetilde{U}_n\vert\mathbf{X}_n}(u)&=&%
\frac{\phi(u,\mathbf{X}_n)}{E[\phi(U_n,\mathbf{X}_n)\vert\mathbf
{X}_n]}f_{U_n\vert\mathbf{X}_n}(u)
,\qquad u>0 .
\end{eqnarray*}
\item In addition,
\begin{eqnarray*}
\phi(\mathbf{X}_n) &=& E\bigl[\phi(U_n,
\mathbf{X}_n)\vert\mathbf{X}_n\bigr] = n+q ,
\\
\omega_{n,\aleph(\mathbf{p}_n)+1}(x)&=&E\bigl[\omega_{n,\aleph(\mathbf
{p}_n)+1}(U_n,
x)\vert\mathbf{X}_n\bigr] \quad \mbox{and}
\\
\omega_{n,k}(Y_k)&=&E\bigl[\omega_{n,k}(U_n,
Y_k)\vert\mathbf{X}_n\bigr], \qquad k=1,\ldots,\aleph(
\mathbf{p}_n) ,%
\end{eqnarray*}
\end{enumerate}
where $\psi_a(b)$ for $a\geq0$ and $b\geq0$ defined in (\ref{psi})
and $\tau_{m,z}(a)$ for $m>0$, $z\in\mathbb{X}$ and $a\geq0$ defined in
(\ref{tau}).
\end{proposition}
\begin{pf}See Appendix \ref{prop1proof}.
\end{pf}
Proposition \ref{prop1BMUrnpredictivedistribution}
gives the predictive distribution and the augmented predictive
distribution for the sequence
$\mathbf{X}_n= \{X_1,\ldots,X_n \}$ as in Blackwell and MacQueen \cite{BM1973} and
James \textit{et al.} \cite{JLP2008}. Statement 1 of Proposition \ref{prop1BMUrnpredictivedistribution} could be viewed as a direct sampling scheme
and statement 2 as a conditional sampling scheme that involves
iterative sampling with $\widetilde{U}_n$ given $\mathbf{X}_n$. From
statement 1, conditional on
$\mathbf{X}_n= \{\mathbf{Y}_{\aleph(\mathbf{p}_n)},\mathbf{p}_n \}$,
$X_{n+1}$ is sampled from
$[\int_{\mathbb{X}}\omega_{n,\aleph(\mathbf{p}_n)+1}(x)\eta
(\mathrm{d}x)]^{-1}\omega_{n,\aleph(\mathbf{p}_n)+1}(x)\eta(\mathrm{d}x)$
with probability
$[\phi(\mathbf{X}_n)]^{-1}\int_{\mathbb{X}}\omega_{n,\aleph(\mathbf
{p}_n)+1}(x)\*\eta(\mathrm{d}x)$
and the new sample is allocated a new index
$\aleph(\mathbf{p}_n)+1,$ becoming $Y_{\aleph(\mathbf{p}_n)+1}$.
Otherwise $X_{n+1}$ has probability $
1-[\phi(\mathbf{X}_n)]^{-1}\int_{\mathbb{X}}\omega_{n,\aleph(\mathbf
{p}_n)+1}(x)\eta(\mathrm{d}x)
$ to be one of the existing sample
$\mathbf{Y}_{\aleph(\mathbf{p}_n)}$. This sequential scheme directly
collects sample of size $n$ cumulatively. Statement 2 of Proposition
\ref{prop1BMUrnpredictivedistribution} suggests an alternative
sampling scheme that draws from
$ \{X_i,\widetilde{U}_i \}$ sequentially for $i=1,\ldots,n$.
To initialize the scheme, a sample of $\widetilde{U}_0$ from the
marginal distribution $f_{\widetilde{U}_0}=f_{U_0}$ and a sample of
$X_1$ from
$[\int_{\mathbb{X}}\tau_{1,x}(\gamma+\widetilde{U}_0)\eta(\mathrm{d}x)]^{-1}\tau_{1,x}(\gamma+\widetilde{U}_0)\eta(\mathrm{d}x)$
are required, then a sample of $ \{X_2,\ldots,X_n \}$ can be
achieved through iterating the following steps
\begin{itemize}
\item Step 1: conditional on $\mathbf{X}_i$, sample $\widetilde{U}_i$ from
$f_{\widetilde{U}_i\vert\mathbf{X}_i}$,
\item Step 2: conditional on $\mathbf{X}_i$ and $\widetilde{U}_i$,
sample $X_{i+1}$ from $\mathbb{P}\{X_{i+1}\in
\mathrm{d}x\vert\mathbf{X}_i,\widetilde{U}_i\}$,
\end{itemize}
for $i=1,\ldots,n-1$, where
$\mathbf{X}_i= \{\mathbf{Y}_{\aleph(\mathbf{p}_i)},\mathbf{p}_i \}$.
The scheme described here is more general than those in existing
articles and provides alternatives to sample from some common
processes, such as the normalized generalized gamma process and
the generalized Dirichlet process (see James \cite{James2002}; Lijoi \textit{et al.} \cite{LMP2005a,LMP2005b} for a direct sampling scheme). These two
cases will be discussed in Section \ref{Sectiontiltedexamples}.
%
\begin{remark}\label{remark3prop1}
The joint distribution of $\mathbf{X}_n$ can be recovered according
to (\ref{tiltedCRMBM}) that can be also found in James \cite{James2002}, Theorem 5.2.
\end{remark}
%
%
\begin{proposition}\label{prop2marginaldistribution} Let
$\widetilde{G}_t=\widetilde\mu_t /\widetilde\mu_t(\mathbb{X})$ be
a normalized tilted \textsc{crm} (Definition \ref{Deftiltedprobabilitymeasure}) defined on $\mathbb{X}$ with $h$ as in
(\ref{generalh}) and $h^\prime$ as in (\ref{PET}), that is
$h(x)=\mathrm{e}^{-\gamma x}x^{-q}$. The parameters of this normalized tilted
\textsc{crm} are given by the measure $\nu$, and the scalers $q\geq0$
and $\gamma\geq0$. Then,
\begin{enumerate}
\item The marginal distribution for
$\mathbf{X}_n= \{\mathbf{Y}_{\aleph(\mathbf{p}_n)},\mathbf{p}_n \}$
is given by
\begin{eqnarray*}
\frac{\int_{\mathbb{R}^{+}}\mathrm{e}^{-\psi_0(\gamma+u)}\prod_{k=1}^{\aleph
(\mathbf{p}_n)}\tau_{n_k,Y_k}(\gamma+u)u^{n+q-1}\,\mathrm{d}u}%
{\sum_{\mathbf{p}_n}\int
_{\mathbb{R}^{+}}\mathrm{e}^{-\psi_0(\gamma+u)}\prod_{k=1}^{\aleph(\mathbf{p}_n)}
\kappa_{n_k}(\gamma+u)u^{n+q-1}\,\mathrm{d}u}%
\prod
_{k=1}^{\aleph(\mathbf{p}_n)}\eta(\mathrm{d}Y_k) ,%
\end{eqnarray*}
where
\begin{eqnarray*}
\kappa_{m}(a)=%
\int_{\mathbb{X}}
\tau_{m,x}(a)\eta(\mathrm{d}x)=%
\int_{\mathbb{X}} \biggl[
\int_{\mathbb{R}^{+}}s^{m}\mathrm{e}^{-sa}%
\rho_x(\mathrm{d}s) \biggr]\eta(\mathrm{d}x) ,\qquad m>0 , a\geq0 .
\end{eqnarray*}
\item The marginal distribution for
$\mathbf{p}_n$ is given by
\begin{eqnarray*}
\frac{\int_{\mathbb{R}^{+}}\mathrm{e}^{-\psi_0(\gamma+u)}\prod_{k=1}^{\aleph
(\mathbf{p}_n)}\kappa_{n_k}(\gamma+u)u^{n+q-1}\,\mathrm{d}u}%
{\sum_{\mathbf{p}_n}
\int_{\mathbb{R}^{+}}\mathrm{e}^{-\psi_0(\gamma+u)}\prod_{k=1}^{\aleph(\mathbf{p}_n)}
\kappa_{n_k}(\gamma+u)u^{n+q-1}\,\mathrm{d}u} .%
\end{eqnarray*}
\item Conditional on $U_n$, the distribution for
$\mathbf{X}_n= \{\mathbf{Y}_{\aleph(\mathbf{p}_n)},\mathbf{p}_n \}$
is given by
\begin{eqnarray*}
\frac{\prod_{k=1}^{\aleph(\mathbf{p}_n)}\tau_{n_k,Y_k}(\gamma+U_n)}%
{\sum_{\mathbf{p}_n}
\prod_{k=1}^{\aleph(\mathbf{p}_n)}\kappa_{n_k}(
\gamma+U_n)}%
\prod_{k=1}^{\aleph(\mathbf{p}_n)}
\eta(\mathrm{d}Y_k) .
\end{eqnarray*}
\item Conditional on $U_n$, the distribution for
$\mathbf{p}_n$ is given by
\begin{eqnarray*}
\frac{\prod_{k=1}^{\aleph(\mathbf{p}_n)}\kappa_{n_k}(\gamma+U_n)}%
{\sum_{\mathbf{p}_n}
\prod_{k=1}^{\aleph(\mathbf{p}_n)}\kappa_{n_k}(
\gamma+U_n)} .%
\end{eqnarray*}
\item Conditional on $\mathbf{p}_n$, the distribution for $U_n$ is
given by
\begin{eqnarray*}
\frac{\mathrm{e}^{-\psi_0(\gamma+u)}\prod_{k=1}^{\aleph(\mathbf{p}_n)}\kappa_{n_k}(\gamma+u)u^{n+q-1}}%
{\int_{\mathbb{R}^{+}}\mathrm{e}^{-\psi_0(\gamma+u)}
\prod_{k=1}^{\aleph(\mathbf
{p}_n)}\kappa_{n_k}(
\gamma+u)u^{n+q-1}\,\mathrm{d}u},%
\end{eqnarray*}
\end{enumerate}
where $\psi_a(b)$ for $a\geq0$ and $b\geq0$ defined in (\ref{psi})
and $\tau_{m,z}(a)$ for $m>0$, $z\in\mathbb{X}$ and $a\geq0$ defined in
(\ref{tau}).
\end{proposition}
\begin{pf}
See Appendix \ref{prop2proof}.
\end{pf}

Proposition \ref{prop2marginaldistribution} gives the
marginal distributions of $\mathbf{X}_n,$ both conditional and
unconditional on $U_n$, in statements 1 and 3, respectively. In
statements 2 and 4, the proposition gives the distributions of
$\mathbf{p}_n$ of first $n$ integers $ \{1,\ldots,n \}$ both
conditional and unconditional on $U_n$. The distributions of the
partitions are the {\it exchangeable partition probability function}
(\textsc{eppf}) as they all are symmetric functions of
$ \{n_1,\ldots,n_{\aleph(\mathbf{p}_n)} \}$. A~special
structure of the \textsc{eppf}, called the Gibbs form, is also
available in the conditional case. Specifically, the distribution of
the partition, $\mathbf{p}_n,$ conditional on $U_n$ (Statement 4 of
Proposition \ref{prop2marginaldistribution}) has the Gibbs
form (Pitman \cite{Pitman2006}, Theorem~4.6, page 86), that is, the \textsc{eppf} is of the form
$V_{n,\aleph(\mathbf{p}_n)}\prod_{k=1}^{\aleph(\mathbf{p}_n)}
W_{n_k}$ in which $V_{n,\aleph(\mathbf{p}_n)}$ is a function of $n$
and $\aleph(\mathbf{p}_n)$ and $W_{n_k}$ is a function of $n_k$.
Such a partition is called the Gibbs partition and therefore
$\mathbf{p}_n$ is a Gibbs partition conditional on $U_n$ (see
Pitman \cite{Pitman1995} and Pitman \cite{Pitman2006} for the details of the
\textsc{eppf} and the Gibbs form).
%
\begin{remark}\label{remark3prop2}
The results appear in Proposition \ref{prop2marginaldistribution}
can be derived from James \cite{James2002}, Theorem~5.2. Some special cases
could be found in James \cite{James2002}, Sections 5.3 and 5.4.
\end{remark}

\section{The tilted version of generalized gamma process and
generalized Dirichlet process}\label{Sectiontiltedexamples}

We consider the tilted version of two interesting
and important classes of measures, namely the generalized gamma
process and the generalized Dirichlet process. Tilting these two
classes of measures yields the normalized beta-gamma process, the
Poisson--Dirichlet process, the normalized generalized gamma process
and the Dirichlet process. They all have been extensively studied
and have wide applications in both statistics and probability. We
now discuss the tilted version of these processes. Let $\widetilde
G_t$ be a normalized tilted \textsc{crm} with $h(x)=\mathrm{e}^{-\gamma x}
x^{-q}$ and with $\nu$ chosen in the following subsections.

\subsection{Generalized gamma process}
The generalized gamma process is
considered a building block for random probability measures
(James \cite{James2002}) and has been widely investigated
(Lijoi \textit{et al.} \cite{LMP2007b}). Earlier studies on this process can be found
in modeling survival functions (Hougaard \cite{Hougaard1986}) and in
spatial modeling (Wolpert and Ickstadt \cite{WI1998};\vadjust{\goodbreak} Brix \cite{Brix1999}). The
intensity of the generalized gamma process is given~by
\begin{eqnarray*}
&&\nu(\mathrm{d}x,\mathrm{d}s)=\eta(\mathrm{d}x)\times\rho_x(\mathrm{d}s) =%
\eta(\mathrm{d}x)\times
\frac{\theta}{\Gamma(1-\alpha)}s^{-1-\alpha}\mathrm{e}^{-b(x)s}\,\mathrm{d}s ,
\\
&&\quad 0\leq\alpha<1 , \theta>0 , b(x)\geq0 .
\end{eqnarray*}
This process is
an important class since it contains the following well known
processes:
\begin{enumerate}
\item Taking $0\leq\alpha<1$, $\theta>0$, $b(x)=0$ yields the intensity
of positive $\alpha$ stable process,
\item Taking $\alpha=0$, $\theta>0$, $b(x)=b\geq0$ yields the intensity
of gamma process, where $b$ is a known constant,
\item Taking $\alpha=0$, $\theta>0$, $b(x)\geq0$ yields the intensity
of extended gamma process and
\item Taking $\alpha=1/2$, $\theta>0$, $b(x)=b\geq0$ yields the
intensity of inverse Gaussian process, where $b$ is a known
constant.
\end{enumerate}
We now consider the normalized random measure derived from the
polynomially and exponentially tilted law of the generalized gamma
process. An application of Theorem \ref{thm2conditionalposteriorG} and Proposition \ref{prop1BMUrnpredictivedistribution} leads
to the conditional \textsc{gbmp} urn formula
%
\begin{eqnarray}
\label{generalgbmp} &&\mathbb{P}\{X_{n+1}\in \mathrm{d}x\vert
\mathbf{X}_{n},\widetilde{U}_n\}\nonumber\\
&&\quad =\biggl(
{\frac{\theta}{(\gamma+\widetilde{U}_n+b(x))^{1-\alpha}}}\nonumber\\
&&\hphantom{\quad =\biggl(}\bigg/
\biggl({\int_{\mathbb{X}}\frac{\theta}{(\gamma+\widetilde{U}_n+b(x))^{1-\alpha
}}
\eta(\mathrm{d}x)+\sum_{k^{\prime}=1}^{\aleph(\mathbf{p}_{n})}
\frac
{n_{k^{\prime}}-\alpha}{\gamma+\widetilde{U}_n+b(Y_{k^{\prime}})}}\biggr)\biggr)
\eta(\mathrm{d}x)
\\
&&\qquad{}+ %
\sum_{k=1}^{\aleph(\mathbf{p}_{n})}%
\biggl({\frac{n_k-\alpha}{\gamma+\widetilde{U}_n+b(Y_k)}}\nonumber\\
&&\hphantom{\qquad{}+ %
\sum_{k=1}^{\aleph(\mathbf{p}_{n})}\biggl(}\bigg/%
\biggl({\int_{\mathbb{X}}
\frac{\theta}{(\gamma+\widetilde{U}_n+b(x))^{1-\alpha
}}\eta(\mathrm{d}x)+\sum_{k^{\prime}=1}^{\aleph(\mathbf{p}_{n})}
\frac
{n_{k^{\prime}}-\alpha}{\gamma+\widetilde{U}_n+b(Y_{k^{\prime}})}}\biggr)\biggr)
\delta_{Y_{k}}(\mathrm{d}x) ,
\nonumber
\end{eqnarray}
and the posterior distribution of $\widetilde{U}_n$ is given by
\begin{eqnarray*}
f_{\widetilde{U}_n\vert\mathbf{X}_n}(u)
\propto\cases{\displaystyle
\Biggl[\int_{\mathbb{X}}
\frac{\theta}{(\gamma+u+b(x))^{1-\alpha}}\eta (\mathrm{d}x)+\sum_{k^{\prime}=1}^{\aleph(\mathbf{p}_{n})}
\frac{n_{k^{\prime
}}-\alpha}{\gamma+u+b(Y_{k^{\prime}})} \Biggr]
\cr
\quad{}\displaystyle\times \frac{\mathrm{e}^{-(\theta/\alpha)\int_{\mathbb{X}}(\gamma+u+b(x))^{\alpha}\eta
(\mathrm{d}x)}u^{n+q}}%
{\prod_{k=1}^{\aleph(\mathbf{p}_{n})}
\bigl(\gamma+u+b(Y_k)\bigr)^{n_k-\alpha}},
\qquad 0<\alpha<1,
\cr
\displaystyle\Biggl[\int_{\mathbb{X}}\frac{\theta}{\gamma+u+b(x)}\eta(\mathrm{d}x)+\sum
_{k^{\prime}=1}^{\aleph(\mathbf{p}_{n})}\frac{n_{k^{\prime}}}{\gamma
+u+b(Y_{k^{\prime}})} \Biggr]
\cr
\quad{}\displaystyle\times \frac{\mathrm{e}^{-\theta\int_{\mathbb{X}}\ln(\gamma+u+b(x))\eta(\mathrm{d}x)}u^{n+q}}%
{\prod_{k=1}^{\aleph(\mathbf{p}_{n})}
\bigl(\gamma+u+b(Y_k)\bigr)^{n_k}},\hspace*{30pt} \alpha=0.
}
\end{eqnarray*}
Notice that when $b(x)=b>0$ and $\alpha=0$, then the condition
(\ref{ConditionEhx})
is reduced to $\theta>q\geq0$. In fact, when $b(x)=b>0$ and
$0<\alpha<1$, then $\theta$ is not required to be greater than $q$ and
the condition (\ref{ConditionEhx}) holds with $q\geq0$ and
$\theta>0$. With a general $b(x)$, the condition (\ref{ConditionEhx}) is required to be examined. By inspection of (\ref{generalgbmp}),
one could find that when $\alpha=0$ and $b(x)=b$,
$X_{n+1}$ is no longer dependent on
$\widetilde{U}_n$ given the past $ \{X_1,\ldots,X_n \}$, that
is $\mathbb{P}\{X_{n+1}\in \mathrm{d}x\vert\mathbf{X}_{n},\widetilde{U}_n\}=\mathbb{P}\{X_{n+1}\in
\mathrm{d}x\vert\mathbf{X}_{n}\}$. This fact is also emphasized in Remark 2 of
James \textit{et al.} \cite{JLP2008}, page 86.
This setting with $q=0$ is corresponding to
the normalized gamma process or the Dirichlet process and
the \textsc{gbmp} urn formula is given by
%
\begin{eqnarray}
\label{bmurn} \mathbb{P}\{X_{n+1}\in \mathrm{d}x\vert\mathbf{X}_{n},
\widetilde{U}_n\}=
\frac{\theta}{\theta+n}
\eta(\mathrm{d}x)+
\sum_{k=1}^{\aleph(\mathbf{p}_{n})}%
\frac{n_k}{\theta+n}
\delta_{Y_{k}}(\mathrm{d}x) .%
\end{eqnarray}
This (\ref{bmurn}) is the \textsc{bmp} urn formula (see also Blackwell and MacQueen \cite{BM1973}).

\subsubsection*{Normalized beta-gamma process}
The law of the beta-gamma
process can be derived from the polynomially tilting the law of the
gamma process (James \cite{James2005b}). The beta-gamma
process was first introduced for the proof of
the Markov--Krein identity of the Dirichlet process mean functionals
and since then it has become a useful analytical tool for studying
the Dirichlet process (see James \cite{James2005b}
and James \textit{et al.} \cite{JLP2008}). Here we mention that Cifarelli and Regazzini
\cite{CR1979a,CR1979b,CR1990} were
the first works in which this identity
was explicitly demonstrated in relation to the research of the law
of the Dirichlet process mean functionals.
An important fact in James \cite{James2005b} states
that a Dirichlet process can be expressed as a normalized
beta-gamma process.
This becomes an interesting alternative expression of the
Dirichlet process that is usually expressed as the normalized gamma process.
This expression is given by taking
$q>0$, $\gamma=0$, $\alpha=0$, $\theta>q$, $b(x)=b\geq0$. The
condition $\theta>q$ (see also James \cite{James2005b}, equation 5, page 649) is equivalent to the condition (\ref{ConditionEhx}). With $\alpha=0$ and $b(x)=b$, as in the construction through the
normalized gamma process, the variable $X_{n+1}$ is not dependent
on $\widetilde{U}_n$ nor $U_n$ given the past
$ \{X_1,\ldots,X_n \}$.

\subsubsection*{Poisson--Dirichlet process} The Poisson--Dirichlet
process is a common and well known process used in both statistical
and probabilistic modeling. This is also called Pitman--Yor process
which is coined by Ishwaran and James \cite{IJ2001}. This process has been shown to
be useful in a variety of interesting applications in combinatorics
(Arratia \textit{et al.} \cite{ABT1997}), population genetics (Griffiths and Lessard \cite{GL2005}) and
Bayesian statistics (Ishwaran and James \cite{IJ2001,IJ2003}). This process was first
introduced by Kingman \cite{Kingman1975} and Pitman and Yor \cite{PY1997} provided a
detailed study of its properties. We consider the Poisson--Dirichlet
process with parameters ($\alpha,q$) that is equivalent to take
$\theta=1$, $\gamma=0$, $b(x)=0$, $q>0$, and $0<\alpha<1$. Special
cases include the Dirichlet process and the normalized stable
process with parameters ($0,q$) and ($\alpha,0$) respectively. These
two processes
could be seen as an two parameter extension of the Dirichlet
process. Even thought there is an explicit expression of the
unconditional \textsc{gbmp}
urn formula (Pitman \cite{Pitman2003,Pitman2006};  Ishwaran and James \cite{IJ2001,IJ2003}),
it is
still worth examining the augmented version for the
Poisson--Dirichlet process ($\alpha,q$) process. The conditional
\textsc{gbmp}
urn formula is given by
\begin{eqnarray*}
\mathbb{P}\{X_{n+1}\in \mathrm{d}x\vert\mathbf{X}_{n},
\widetilde{U}_n\}=%
\frac{\widetilde{U}_n^{\alpha}}{\widetilde{U}_n^{\alpha}+n-\alpha\aleph
(\mathbf{p}_{n})}%
\eta(\mathrm{d}x)+%
\sum_{k=1}^{\aleph(\mathbf{p}_{n})}
\frac{n_{k}-\alpha}{ \widetilde{U}_n^{\alpha}+n-\alpha\aleph(\mathbf
{p}_{n})}%
\delta_{Y_{k}}(\mathrm{d}x) ,
\end{eqnarray*}
where
%
\begin{eqnarray}
\label{UgivenXPD} f_{\widetilde{U}_n\vert\mathbf{X}_n}(u)&=& \frac{q+\alpha\aleph(\mathbf{p}_{n})}%
{n+q}%
\frac{\mathrm{e}^{-(1/\alpha)u^{\alpha}}u^{q+\alpha\aleph(\mathbf{p}_{n})+\alpha
-1}}%
{\alpha^{q/\alpha+\aleph(\mathbf{p}_{n})}\Gamma
(q/\alpha+\aleph(\mathbf {p}_{n})+1)}\nonumber
\\[-8pt]\\[-8pt]
&&{}+ %
\frac{n-\alpha\aleph(\mathbf{p}_{n})}%
{n+q}%
\frac{\mathrm{e}^{-(1/\alpha)u^{\alpha}}u^{q+\alpha\aleph(\mathbf{p}_{n})-1}}%
{\alpha^{q/\alpha+\aleph(\mathbf{p}_{n})-1}\Gamma(q/\alpha+\aleph (
\mathbf{p}_{n}))} .
\nonumber
\end{eqnarray}
This is equivalent to $\widetilde{U}_n\overset{d}{=}G_1^{1/\alpha}$
with probability $(q+\alpha\aleph(\mathbf{p}_{n}))/(n+q)$ and
$\widetilde{U}_n\overset{d}{=}G_2^{1/\alpha}$ with probability
$(n-\alpha\aleph(\mathbf{p}_{n}))/(n+q)$ where $G_1$ is a
Gamma($q/\alpha+\aleph(\mathbf{p}_{n})+1,1/\alpha$) random variable
and $G_2$ is a Gamma($q/\alpha+\aleph(\mathbf{p}_{n}),1/\alpha$)
random variable.

\subsubsection*{Normalized generalized gamma process} The
normalized generalized gamma process with $\theta>0$, $\gamma=0$,
$b(x)=b>0$, $q=0$, and $0<\alpha<1$ is considered in
James \cite{James2002} and Lijoi \textit{et al.} \cite{LMP2005b,LMP2007b}. Specifically, the
conditional \textsc{gbmp}
urn formula is given by
\begin{eqnarray*}
\mathbb{P}\{X_{n+1}\in \mathrm{d}x\vert\mathbf{X}_{n},
\widetilde{U}_n\} &=&
 \frac{\theta(\widetilde{U}_n+b)^{\alpha}}{\theta(\widetilde
{U}_n+b)^{\alpha}+n-\alpha\aleph(\mathbf{p}_{n})}%
\eta(\mathrm{d}x)\\
&&{}+ \sum_{k=1}^{\aleph(\mathbf{p}_{n})}
\frac{n_{k}-\alpha}{\theta(\widetilde{U}_n+b)^{\alpha}+n-\alpha\aleph
(\mathbf{p}_{n})}%
\delta_{Y_{k}}(\mathrm{d}x) ,
\end{eqnarray*}
where the density of the augmented variable $\widetilde{U}_n$ is
given by
%
\begin{eqnarray}
\label{UgivenXNGG} f_{\widetilde{U}_n\vert\mathbf{X}_n}(u)%
\propto \frac{%
[\theta(u+b)^{\alpha}+n-\alpha\aleph(\mathbf{p}_{n})]%
\mathrm{e}^{-(\theta/\alpha)(u+b)^{\alpha}} u^{n}%
}{%
(u+b)^{n-\alpha\aleph(\mathbf{p}_{n})+1}%
}
.%
\end{eqnarray}
A little effort might be required for sampling $\widetilde{U}_n$ from
(\ref{UgivenXNGG}).
One could follow Devroye \cite{Devroye1986}, Section II.3.3, page 47, to derive a suitable rejection procedure. Here we
give a simple illustration. A sample of $V_n$ could be drawn from
\begin{eqnarray*}
\frac{\alpha
\aleph(\mathbf{p}_{n})}{n}\frac{\theta^{\aleph(\mathbf{p}_{n})+1}%
\mathrm{e}^{-(\theta/\alpha)v^{\alpha}}v^{\alpha\aleph(\mathbf{p}_{n})+\alpha-1}%
}{\alpha^{\aleph(\mathbf{p}_{n})}\Gamma(\aleph(\mathbf{p}_{n})+1)}%
+\frac{n-\alpha\aleph(\mathbf{p}_{n})}{n}
\frac{\theta^{\aleph(\mathbf
{p}_{n}%
)}\mathrm{e}^{-(\theta/\alpha) v^{\alpha}}v^{\alpha\aleph(\mathbf{p}_{n})-1}}%
{\alpha^{\aleph(\mathbf{p}_{n})-1}\Gamma(\aleph(
\mathbf{p}_{n}))} ,%
\end{eqnarray*}
and if $\zeta< \psi(V_n)$, then $\widetilde{U}_n=V_n$, otherwise
sample $V_n$ again until $\zeta< \psi(V_n)$ where $\zeta$ is an
uniform random variable which is independent of $V_n$ and $\psi(v)=
\mathrm{e}^{-(\theta/\alpha)[(v+b)^{\alpha}-v^{\alpha}]}\times(\frac{v}%
{v+b})^{n-\alpha\aleph(\mathbf{p}_{n})+1}\times\frac{\theta
(v+b)^{\alpha}%
+(n-\alpha\aleph(\mathbf{p}_{n}))}{\theta v^{\alpha}+(n-\alpha\aleph
(\mathbf{p}_{n}))}$. Notice that
$\widetilde U_n$ given $\mathbf{X}_n$ (\ref{UgivenXNGG}) and $V_n$ are identical in distribution when $b=0$.
The random variable $V_n$ can be described as:
$V_n\overset{d}{=}G_1^{1/\alpha}$ with probability
$\alpha\aleph(\mathbf{p}_{n})/n$ and
$V_n\overset{d}{=}G_2^{1/\alpha}$ with probability
$(n-\alpha\aleph(\mathbf{p}_{n}))/n$ where $G_1$ is a
Gamma($\aleph(\mathbf{p}_{n})+1,\theta/\alpha$) random variable and
$G_2$ is a $\Gam$($\aleph(\mathbf{p}_{n}),\theta/\alpha$) random
variable.

\subsection{Generalized Dirichlet process} Regazzini \textit{et al.} \cite{RLP2003} introduce
the generalized Dirichlet process as an example
for determining the mean of normalized random measures with
independent increments. Apart from studying probabilistic properties
of the generalized Dirichlet process, its use in Bayesian
non-parametric statistics is developed in Lijoi \textit{et al.} \cite{LMP2005a}. We state
the intensity of the generalized Dirichlet process as
\begin{eqnarray*}
\nu(\mathrm{d}s,\mathrm{d}x)=\eta(\mathrm{d}x) \times\rho_x(\mathrm{d}s) =\eta(\mathrm{d}x) \times \theta
\frac{1-\mathrm{e}^{-cs}}{1-\mathrm{e}^{-s}}s^{-1}\mathrm{e}^{-s}\,\mathrm{d}s ,
\end{eqnarray*}
and introduce the difference of two Hurwitz Zeta functions as,
\begin{eqnarray*}
\varphi_{n_k}(\gamma+u,c)%
=\sum
_{\ell=0}^{\infty} \biggl[\frac{1}{(\gamma+u+1+\ell)^{n_{k}}}-
\frac
{1}{(\gamma+u+c+1+\ell)^{n_{k}}} \biggr] .
\end{eqnarray*}
When $c$ is a positive integer (which is considered by
Regazzini \textit{et al.} \cite{RLP2003} and Lijoi \textit{et al.} \cite{LMP2005a}), the function can be
simplified to a finite sum,
\begin{eqnarray*}
\varphi_{n_k}(\gamma+u,c)%
=\sum
_{\ell=0}^{c-1}\frac{1}{(\gamma+u+1+\ell)^{n_{k}}} ,\qquad c=1,2,\ldots.
\end{eqnarray*}
The conditional \textsc{gbmp}
urn is given by
\begin{eqnarray*}
\mathbb{P}\{X_{n+1}\in \mathrm{d}x\vert\mathbf{X}_{n},
\widetilde{U}_n\}&=&%
\biggl({\theta\varphi_{1}(\gamma+\widetilde{U}_{n},c)}\\
&&\hspace*{5pt}\bigg/\biggl({\theta
\varphi_{1}(\gamma+\widetilde{U}_{n},c)+%
\sum
_{k^{\prime}=1}^{\aleph(\mathbf{p}_{n})}n_{k^{\prime}}%
\frac{\varphi_{n_{k^{\prime}}+1}(\gamma+\widetilde{U}_{n},c)}{\varphi_{n_{k^{\prime}}}(\gamma+\widetilde{U}_{n},c)}}\biggr)\!\biggr)%
\eta(\mathrm{d}x)%
\\
&&{} + %
\sum_{k=1}^{\aleph(\mathbf{p}_{n})}%
\biggl({n_{k}\frac{\varphi_{n_{k}+1}(\gamma+\widetilde{U}_{n},c)}{\varphi_{n_{k}}(\gamma+\widetilde{U}_{n},c)}}\\
&&\hspace*{35pt}\bigg/\biggl({\theta\varphi_{1}(\gamma+
\widetilde{U}_{n},c)+%
\sum_{k^{\prime}=1}^{\aleph(\mathbf{p}_{n})}n_{k^{\prime}}%
\frac{\varphi_{n_{k^{\prime}}+1}(\gamma+\widetilde{U}_{n},c)}{\varphi_{n_{k^{\prime}}}(\gamma+\widetilde{U}_{n},c)}}\biggr)\!\biggr)%
\delta_{Y_{k}}(\mathrm{d}x) ,
\end{eqnarray*}
where the density of the augmented variable is given by
\begin{eqnarray*}
f_{\widetilde{U}_{n}\vert\mathbf{X}_{n}}(u) &\propto& %
\chi(u)\\ &=& \Biggl[ \theta
\varphi_{1}(\gamma+u,c)+%
\sum_{k^{\prime}=1}^{\aleph(\mathbf{p}_{n})}n_{k^{\prime}}%
\frac{\varphi_{n_{k^{\prime}}+1}(\gamma+u,c)}{\varphi_{n_{k^{\prime
}}}(\gamma+u,c)} \Biggr]
\\
&&\times \biggl[\frac{\Gamma(\gamma+u+1)}{\Gamma(\gamma+u+c+1)} \biggr]^\theta%
u^{n+q}
\times%
\theta^{\aleph(\mathbf{p}_{n})}%
\prod
_{k=1}^{\aleph(\mathbf{p}_{n})}%
\varphi_{n_k}(\gamma+u,c) .
\end{eqnarray*}
Again, the condition (\ref{ConditionEhx}) stated in the
construction that
$E [h (\mu(\mathbb{X}) ) ]<\infty$ is equivalent
to $\theta>q$. In particular, when $c=1$ this is corresponding to the
Dirichlet process,
$\varphi_{1}(\gamma+u,1)=\varphi_{n_{k}+1}(\gamma+u,1)  /
\varphi_{n_{k}}(\gamma+u,1)$ and it can be shown that $X_{n+1}$
given $ \{X_1,\ldots,X_n \}$ is not dependent on $U_n$. Similar
to the normalized generalized gamma process,
the rejection method can also be proposed due to Devroye \cite{Devroye1986}, Section II.3.3, page 47. That is, sample $W_n$ from
$\Bet(n+q+1,\theta-q)$, and if $\zeta<
\psi((\gamma+1)W_n/(1-W_n))$, then
$\widetilde{U}_n=(\gamma+1)W_n/(1-W_n)$, otherwise sample $W_n$
again until $\zeta< \psi((\gamma+1)W_n/(1-W_n))$ where $\zeta$ is
an uniform random variable which is independent of $W_n$ and
$\psi(v)=\chi(v)\times(\gamma+v+1)^{\theta+n+1}/(c^{\aleph(\mathbf
{p}_{n})}(c\theta+n))$.

\section{Simulation study on approximating distribution of partition
size via the augmented Blackwell--MacQueen P\'olya urn formula}\label{Sectionsim}

We conduct a simulation study on approximating the posterior
probabilities of partition sizes,
$\mathbb{P} \{\aleph(\mathbf{p}_n)=i \}$ for $i=1,\ldots,n$,
using the conditional \textsc{gbmp}
urn formula
discussed in last section. We consider two popular random
probability measures, the Poisson--Dirichlet process and the
normalized generalized gamma process. For the Poisson--Dirichlet
process, we simulate data according to $\alpha=0.5$, $\theta=1$,
$q=1$, $\gamma=0$, $b(x)=0.$ For the normalized generalised gamma
process, we simulate data according to $\alpha=0.5$, $\theta=1$,
$q=0$, $\gamma=0$, $b(x)=1.$ In both cases, we set $n=50.$ We
examine two exact sequential sampling schemes and two MCMC schemes.
Specifically:

\begin{enumerate}[A.2.]
\item[A.1.] Sample $X_i$ sequentially for $i=1,\ldots,n$ according to the
unconditional \textsc{gbmp}
urn formula,
$\mathbb{P} \{X_{i+1}\in \mathrm{d}x\vert\mathbf{X}_{i}  \}$ for
$i=1,\ldots, n-1$.
\item[A.2.] Sample $ \{X_i,\widetilde{U}_i \}$ sequentially for
$i=1,\ldots,n$ according to the
conditional \textsc{gbmp}
urn formula,
$\mathbb{P} \{X_{i+1}\in \mathrm{d}x\vert\mathbf{X}_{i},\widetilde{U}_i \}$ and
$f_{\widetilde{U}_i\vert\mathbf{X}_{i}}(u)$ for $i=1,\ldots, n-1$.
\item[A.3.] Re-sample $X_i$ iteratively for $i=1,\ldots,n$ according
to the
unconditional \textsc{gbmp}
urn formula,
$\mathbb{P} \{X_{i}\in \mathrm{d}x\vert\mathbf{X}_{n}\backslash\{X_i\}
\}$ for $i=1,\ldots, n$.
\item[A.4.] Re-sample $ \{X_i,U_n \}$ iteratively for $i=1,\ldots
,n$ according to the
conditional \textsc{gbmp}
urn formula,
$\mathbb{P} \{X_{i}\in \mathrm{d}x\vert
\mathbf{X}_{n}\backslash\{X_i\},U_n  \}$ and
$f_{U_n\vert\mathbf{X}_{n}}(u)$ for $i=1,\ldots, n$.
\end{enumerate}
Algorithms A.1 and A.2 are exact and are identical in
distribution. Algorithm A.1 has been frequently used in the
literature and for the normalized generalized\vadjust{\goodbreak} gamma process requires
the evaluation of the complicated functions.
Algorithm A.2 is the conditional
\textsc{gbmp}
urn formula derived from the tilted
measure proposed in this article. This method is straightforward to
implement without
much complication in evaluations.
We also include two MCMC Gibbs sampling algorithms described in
A.3 and A.4 for comparison. These are not exact sampling
algorithms and an initial sampling period is necessary to converge
to the stationary distribution. The stationary distribution itself
is identical to that of A.1 and A.2.

\begin{figure}

\includegraphics{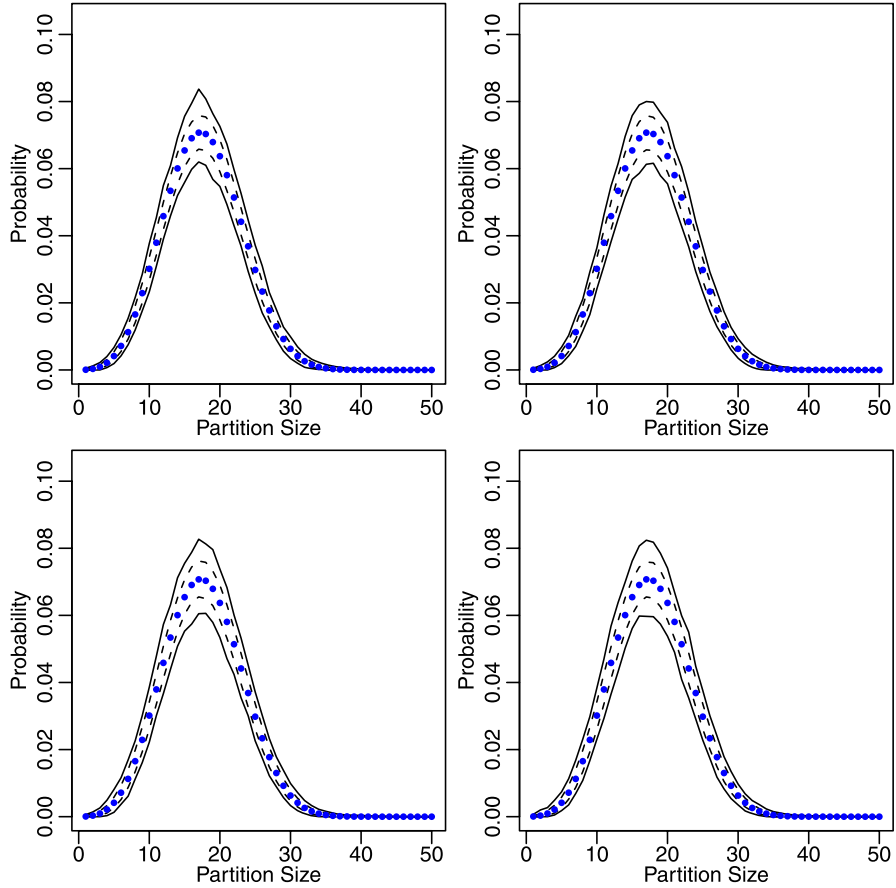}

\caption{Approximate probabilities of partition sizes of the
Poisson Dirichlet process ($\alpha=0.5$, $\theta=1$, $q=1$,
$\gamma=0$, $b(x)=0$),
$\mathbb{P} \{\aleph(\mathbf{p}_{50})=i \}$ for
$i=1,\ldots,50$.
Top Left: A.1 algorithm; %
Top Right: A.2 algorithm; %
Bottom Left: A.3 algorithm; %
Bottom Right: A.4 algorithm. %
The solid bound lines indicate the range of all 10,000 approximates
of the probabilities; The dash bound lines indicate
the 95\% confidence level (2.5\% and 97.5\% quantiles) of
all 10,000 approximates of the probabilities. The dots indicate the
true probabilities.}\label{fig1}
\end{figure}

\begin{figure}

\includegraphics{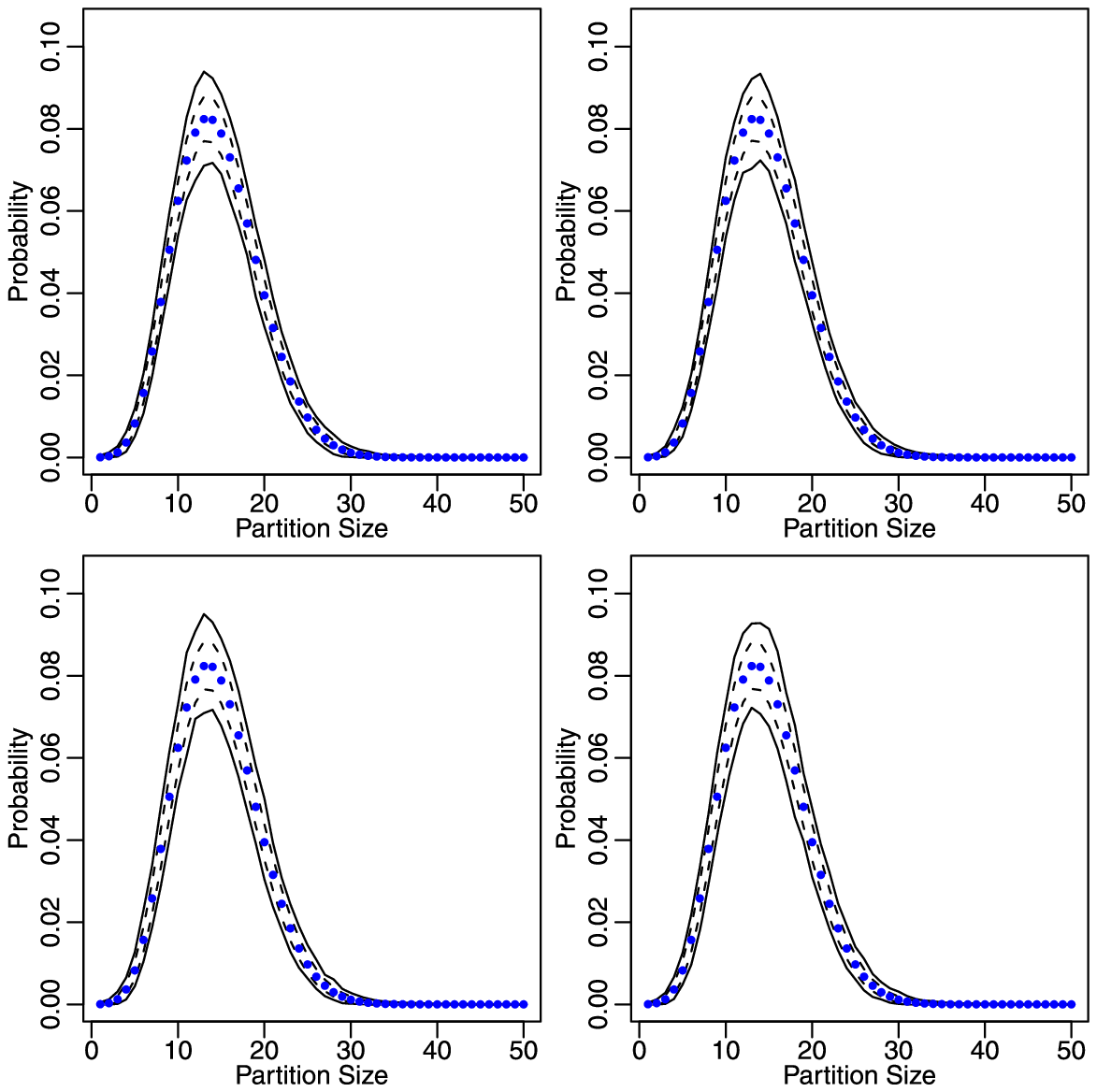}

\caption{Approximate probabilities of partition sizes of the
normalized generalized Gamma process ($\alpha=0.5$, $\theta=1$,
$q=0$, $\gamma=0$, $b(x)=1$),
$\mathbb{P} \{\aleph(\mathbf{p}_{50})=i \}$ for
$i=1,\ldots,50$. %
Top Left: A.1 algorithm; %
Top Right: A.2 algorithm; %
Bottom Left: A.3 algorithm; %
Bottom Right: A.4 algorithm. %
The solid bound lines indicate the range of all 10,000 approximates
of the probabilities; The dash bound lines indicate
the 95\% confidence level (2.5\% and 97.5\% quantiles) of
all 10,000 approximates of the probabilities. The dots indicate the
true probabilities.}\label{fig2}
\end{figure}

In each replication, we sample $10\mbox{,}000$
independent samples from A.1 and A.2 to approximate
$\mathbb{P} \{\aleph(\mathbf{p}_n)=i \}$ for $i=1,\ldots,n$.
Starting with a partition
$\mathbf{p}_n= \{\{1\},\ldots,\{n\} \}$ with all singleton
clusters, we draw $20\mbox{,}000$ samples from algorithms A.3 and A.4.
We ignore the first $10\mbox{,}000$ warmup samples and use the last
$10\mbox{,}000$ samples to approximate
$\mathbb{P} \{\aleph(\mathbf{p}_n)=i \}$ for $i=1,\ldots,n$.
So, each
algorithm produces $10\mbox{,}000$ approximates of probabilities of
partition sizes,
$\mathbb{P} \{\aleph(\mathbf{p}_n)=1 \},\ldots,\mathbb{P} \{
\aleph(\mathbf{p}_n)=n \}$.
To summarize the results, Figures \ref{fig1} and \ref{fig2} show the range and
95\% confidence levels of $10\mbox{,}000$ approximates of probabilities of
partition sizes for the Poisson--Dirichlet process and the normalized
generalized gamma process respectively for algorithms A.1--A.4.
Similarly, Tables \ref{tab1} and \ref{tab2} shows the true probabilities and the
means and standard errors of the approximates given by algorithms
A.1--A.4.

%

Figure \ref{fig1} and Table \ref{tab1} show that for the Poisson--Dirichlet process
algorithms A.1 and A.2 result in samples from identical
distributions, as the theory would suggest. The MCMC algorithms,
A.3 and A.4 also produce similar results to A.1 and A.2
except for the standard errors in Table \ref{tab1}. The standard errors
indicate the variability of the MCMC generated samples is generally
greater than the exact sequential sampling algorithms, as we would
expect. The results for the normalized generalized gamma process
shown in Figure \ref{fig2} and Table \ref{tab2} point to a similar story as in the
Poisson--Dirichlet process.

Finally, we note that it is not necessary to sample $\mathbf{X}_n$
to conduct the simulation. It could be done by simply simulating
partitions $\mathbf{p}_n$ directly instead. It is possible to
integrate out all $\mathbf{Y}_{\aleph(\mathbf{p}_n)}$ from the
conditional \textsc{gbmp}
urn formula and obtain the
weights for partition sampling using the Chinese restaurant process
(see Aldous \cite{Aldous1985}, Lo \textit{et al.} \cite{LBC1996} and
Pitman \cite{Pitman2006}).

\section{Conclusion and further research}
\label{SectionConclusionRemarks}

This article has introduced a class of random probability measures
based on polynomially and exponentially tilting. We have provided a
complete Bayesian analysis of this class of measures with details on
the prior and posterior laws and shown that the class is
structurally conjugate. We described a conditional
Blackwell--MacQueen P\'olya urn sampling scheme that simplifies the
computational requirements to implement such sampling schemes. The
new sampling scheme yields similar answers to more complicated
schemes described in the literature.

We also note a general tilting treatment could be considered for any
measurable function $h$ on $\mathbb{R}^+$ evaluated at the total
mass of the corresponding \textsc{crm} in (\ref{tiltedlawsacledh}).
This general class of random probability measures with homogeneous
intensity, $\nu(\mathrm{d}x,\mathrm{d}s)=\eta(\mathrm{d}x)\rho(\mathrm{d}s)$, can be shown to be the
Poisson--Kingman process (Pitman \cite{Pitman2003,Pitman2006}) by showing
that the normalized random probability measure (\ref{tiltedlawsacledh}) has a Poisson--Kingman partition ($\rho,\gamma$),
where $\gamma$ denotes the tilted density of total mass such that
$\gamma=f_{\widetilde{T}_t}$ (\ref{densityT})
and $\rho$
represents the $\sigma$-finite non-atomic intensity measure $\rho(\mathrm{d}s)$. Then,
the conditional partition distribution is given by
%
\begin{eqnarray}
\label{partitiondistributiongivent} \mathbb{P}(\mathbf{p}_n\vert t)=
\int_{(\mathbb{R}^+)^{\aleph(\mathbf{p}_n)}}%
\mathbb{I}_{\{{t}-\sum_{k=1}^{\aleph(\mathbf{p}_n)}s_k>0\}}
\frac{%
f_{\widetilde{T}} ({t}-\sum_{k=1}^{\aleph(\mathbf{p}_n)}s_k )%
}%
{%
{t}^{n}f_{\widetilde{T}}({t})
}%
\prod_{k=1}^{\aleph(\mathbf{p}_n)}s_k^{n_k}
\rho(\mathrm{d}s_k) .%
\end{eqnarray}
The proof is given in Appendix \ref{partitiondistributiongiventproof}; See also Pitman \cite{Pitman2003}, Lemma 5, page 8, for the related
joint distribution of partitions and jumps. Investigation of the
properties of this general class of random probably measures and
even special cases of (\ref{tiltedlawsacledh}) is interesting for
further research.

%
\begin{sidewaystable*}
\tablewidth=\textwidth
\tabcolsep=0pt
\caption{{Approximate probabilities of partition sizes of the
Poisson Dirichlet process ($\alpha=0.5$, $\theta=1$, $q=1$,
$\gamma=0$, $b(x)=0$),
$\mathbb{P} \{\aleph(\mathbf{p}_{50})=i \}$ for
$i=1,\ldots,50$. %
\emph{First Column}: Partition sizes; %
\emph{Second Column}: True probability of partition sizes; %
\emph{Third Column}: Mean of 10,000 approximates of the probabilities
according to algorithm A.1; %
\emph{Fourth Column}: Standard error of 10,000 approximates of the
probabilities according to algorithm A.1; %
\emph{Fifth Column}: Mean of 10,000 approximates of the probabilities
according to algorithm A.2; %
\emph{Sixth Column}: Standard error of 10,000 approximates of the
probabilities according to algorithm A.2; %
\emph{Seventh Column}: Mean of 10,000 approximates of the probabilities
according to algorithm A.3; %
\emph{Eighth Column}: Standard error of 10,000 approximates of the
probabilities according to algorithm A.3; %
\emph{Ninth Column}: Mean of 10,000 approximates of the probabilities
according to algorithm A.4; %
\emph{Tenth Column}: Standard error of 10,000 approximates of the
probabilities according to algorithm A.4}}\label{tab1}
\begin{tabular*}{\textwidth}{@{\extracolsep{\fill}}llllllllll@{}}
\hline
& $\mathbb{P} \{\aleph(\mathbf{p}_n)=i \}$ & \multicolumn{2}{l}{{A.1 algorithm}}
& \multicolumn{2}{l}{ {A.2 algorithm}} & \multicolumn{2}{l}{{A.3 algorithm}} &
\multicolumn{2}{l}{{A.4 algorithm}} \\[-5pt]
&& \multicolumn{2}{l}{\hrulefill} & \multicolumn{2}{l}{\hrulefill} &
\multicolumn{2}{l}{\hrulefill} & \multicolumn{2}{l}{\hrulefill}\\
& {True} &  {Mean of the} &  {SE of the} &  {Mean of the} &
{SE of the} &  {Mean of the} &  {SE of the} &  {Mean of the} &  {SE of the} \\
$i$ &  {probability} &  {approximates} &  {approximates} & {approximates} &  {approximates}
&  {approximates} &  {approximates} &  {approximates} &  {approximates} \\
\hline
\hphantom{0}1 & 0.000063 & 0.000062 & 0.000080 & 0.000063 & 0.000079 & 0.000063 &
0.000088 & 0.000064 & 0.000088 \\
\hphantom{0}2 & 0.000315 & 0.000316 & 0.000179 & 0.000318 & 0.000179 & 0.000313 &
0.000207 & 0.000311 & 0.000207 \\
\hphantom{0}3 & 0.000936 & 0.000933 & 0.000306 & 0.000934 & 0.000308 & 0.000934 &
0.000369 & 0.000936 & 0.000371 \\
\hphantom{0}4 & 0.002139 & 0.002138 & 0.000462 & 0.002133 & 0.000460 & 0.002138 &
0.000567 & 0.002141 & 0.000575 \\
\hphantom{0}5 & 0.004142 & 0.004146 & 0.000639 & 0.004140 & 0.000641 & 0.004147 &
0.000806 & 0.004135 & 0.000807 \\
\hphantom{0}6 & 0.007139 & 0.007129 & 0.000837 & 0.007143 & 0.000836 & 0.007139 &
0.001071 & 0.007131 & 0.001068 \\
\hphantom{0}7 & 0.011258 & 0.011256 & 0.001060 & 0.011242 & 0.001066 & 0.011257 &
0.001351 & 0.011258 & 0.001333 \\
\hphantom{0}8 & 0.016537 & 0.016519 & 0.001274 & 0.016563 & 0.001277 & 0.016536 &
0.001613 & 0.016526 & 0.001610 \\
\hphantom{0}9 & 0.022898 & 0.022916 & 0.001510 & 0.022901 & 0.001502 & 0.022883 &
0.001868 & 0.022890 & 0.001849 \\
10 & 0.030135 & 0.030133 & 0.001705 & 0.030123 & 0.001716 & 0.030089 &
0.002075 & 0.030119 & 0.002101 \\
11 & 0.037923 & 0.037944 & 0.001923 & 0.037937 & 0.001905 & 0.037893 &
0.002248 & 0.037910 & 0.002311 \\
12 & 0.045836 & 0.045829 & 0.002081 & 0.045820 & 0.002121 & 0.045812 &
0.002430 & 0.045901 & 0.002421 \\
13 & 0.053388 & 0.053358 & 0.002248 & 0.053434 & 0.002243 & 0.053396 &
0.002515 & 0.053364 & 0.002539 \\
14 & 0.060073 & 0.060076 & 0.002407 & 0.060042 & 0.002371 & 0.060038 &
0.002572 & 0.060061 & 0.002600 \\
15 & 0.065424 & 0.065423 & 0.002486 & 0.065396 & 0.002497 & 0.065402 &
0.002638 & 0.065415 & 0.002683 \\
16 & 0.069059 & 0.069067 & 0.002527 & 0.069071 & 0.002554 & 0.069094 &
0.002677 & 0.069051 & 0.002689 \\
17 & 0.070723 & 0.070738 & 0.002575 & 0.070720 & 0.002555 & 0.070754 &
0.002673 & 0.070746 & 0.002684 \\
18 & 0.070317 & 0.070311 & 0.002544 & 0.070328 & 0.002558 & 0.070302 &
0.002729 & 0.070290 & 0.002691 \\
19 & 0.067906 & 0.067983 & 0.002499 & 0.067911 & 0.002480 & 0.067937 &
0.002692 & 0.067889 & 0.002688 \\
20 & 0.063706 & 0.063723 & 0.002443 & 0.063674 & 0.002479 & 0.063756 &
0.002640 & 0.063700 & 0.002636 \\
21 & 0.058061 & 0.058033 & 0.002325 & 0.058056 & 0.002349 & 0.058063 &
0.002585 & 0.058066 & 0.002583 \\
\hline
\end{tabular*}
\end{sidewaystable*}
\setcounter{table}{0}

\begin{sidewaystable*}
\tablewidth=\textwidth
\tabcolsep=0pt
\caption{(Continued)}
%
\begin{tabular*}{\textwidth}{@{\extracolsep{\fill}}llllllllll@{}}
\hline
& $\mathbb{P} \{\aleph(\mathbf{p}_n)=i \}$ & \multicolumn{2}{l}{{A.1 algorithm}}
& \multicolumn{2}{l}{ {A.2 algorithm}} & \multicolumn{2}{l}{{A.3 algorithm}} &
\multicolumn{2}{l}{{A.4 algorithm}} \\[-5pt]
&& \multicolumn{2}{l}{\hrulefill} & \multicolumn{2}{l}{\hrulefill} &
\multicolumn{2}{l}{\hrulefill} & \multicolumn{2}{l}{\hrulefill}\\
& {True} &  {Mean of the} &  {SE of the} &  {Mean of the} &
{SE of the} &  {Mean of the} &  {SE of the} &  {Mean of the} &  {SE of the} \\
$i$ &  {probability} &  {approximates} &  {approximates} & {approximates} &  {approximates}
&  {approximates} &  {approximates} &  {approximates} &  {approximates} \\
\hline
22 & 0.051397 & 0.051395 & 0.002197 & 0.051413 & 0.002207 & 0.051427 &
0.002442 & 0.051390 & 0.002478 \\
23 & 0.044176 & 0.044171 & 0.002073 & 0.044176 & 0.002040 & 0.044161 &
0.002278 & 0.044206 & 0.002309 \\
24 & 0.036847 & 0.036842 & 0.001887 & 0.036834 & 0.001885 & 0.036833 &
0.002103 & 0.036846 & 0.002140 \\
25 & 0.029805 & 0.029810 & 0.001705 & 0.029819 & 0.001708 & 0.029830 &
0.001912 & 0.029824 & 0.001945 \\
26 & 0.023361 & 0.023352 & 0.001510 & 0.023346 & 0.001515 & 0.023366 &
0.001707 & 0.023369 & 0.001712 \\
27 & 0.017724 & 0.017689 & 0.001306 & 0.017739 & 0.001320 & 0.017739 &
0.001477 & 0.017737 & 0.001493 \\
28 & 0.013001 & 0.012989 & 0.001127 & 0.013005 & 0.001131 & 0.013004 &
0.001266 & 0.013002 & 0.001259 \\
29 & 0.009208 & 0.009204 & 0.000958 & 0.009221 & 0.000963 & 0.009200 &
0.001042 & 0.009208 & 0.001072 \\
30 & 0.006287 & 0.006294 & 0.000783 & 0.006280 & 0.000801 & 0.006282 &
0.000859 & 0.006291 & 0.000858 \\
31 & 0.004130 & 0.004141 & 0.000645 & 0.004137 & 0.000639 & 0.004129 &
0.000688 & 0.004138 & 0.000690 \\
32 & 0.002606 & 0.002606 & 0.000513 & 0.002608 & 0.000514 & 0.002606 &
0.000537 & 0.002603 & 0.000539 \\
33 & 0.001575 & 0.001574 & 0.000395 & 0.001571 & 0.000398 & 0.001579 &
0.000423 & 0.001576 & 0.000409 \\
34 & 0.000910 & 0.000906 & 0.000302 & 0.000908 & 0.000299 & 0.000911 &
0.000308 & 0.000913 & 0.000312 \\
35 & 0.000501 & 0.000502 & 0.000225 & 0.000499 & 0.000223 & 0.000498 &
0.000227 & 0.000502 & 0.000228 \\
36 & 0.000261 & 0.000261 & 0.000162 & 0.000261 & 0.000160 & 0.000258 &
0.000163 & 0.000263 & 0.000164 \\
37 & 0.000129 & 0.000129 & 0.000113 & 0.000130 & 0.000113 & 0.000128 &
0.000113 & 0.000129 & 0.000114 \\
38 & 0.000060 & 0.000061 & 0.000078 & 0.000061 & 0.000079 & 0.000060 &
0.000078 & 0.000059 & 0.000078 \\
39 & 0.000026 & 0.000026 & 0.000051 & 0.000027 & 0.000052 & 0.000027 &
0.000052 & 0.000027 & 0.000051 \\
40 & 0.000011 & 0.000011 & 0.000032 & 0.000010 & 0.000032 & 0.000011 &
0.000033 & 0.000010 & 0.000032 \\
41 & 0.000004 & 0.000004 & 0.000020 & 0.000004 & 0.000020 & 0.000003 &
0.000019 & 0.000004 & 0.000021 \\
42 & 0.000001 & 0.000001 & 0.000012 & 0.000001 & 0.000012 & 0.000001 &
0.000012 & 0.000001 & 0.000012 \\
43 & 0.000000 & 0.000000 & 0.000006 & 0.000000 & 0.000006 & 0.000000 &
0.000006 & 0.000000 & 0.000006 \\
44 & 0.000000 & 0.000000 & 0.000004 & 0.000000 & 0.000004 & 0.000000 &
0.000003 & 0.000000 & 0.000003 \\
45 & 0.000000 & 0.000000 & 0.000002 & 0.000000 & 0.000002 & 0.000000 &
0.000002 & 0.000000 & 0.000002 \\
46 & 0.000000 & 0.000000 & 0.000000 & 0.000000 & 0.000000 & 0.000000 &
0.000000 & 0.000000 & 0.000000 \\
47 & 0.000000 & 0.000000 & 0.000000 & 0.000000 & 0.000000 & 0.000000 &
0.000000 & 0.000000 & 0.000000 \\
48 & 0.000000 & 0.000000 & 0.000000 & 0.000000 & 0.000000 & 0.000000 &
0.000000 & 0.000000 & 0.000000 \\
49 & 0.000000 & 0.000000 & 0.000000 & 0.000000 & 0.000000 & 0.000000 &
0.000000 & 0.000000 & 0.000000 \\
50 & 0.000000 & 0.000000 & 0.000000 & 0.000000 & 0.000000 & 0.000000 &
0.000000 & 0.000000 & 0.000000 \\
\hline
\end{tabular*}
\end{sidewaystable*}

\begin{sidewaystable*}
\tablewidth=\textwidth
\tabcolsep=0pt
\caption{Approximate of probabilities of partition sizes of the
normalized generalized Gamma process ($\alpha=0.5$, $\theta=1$,
$q=0$, $\gamma=0$, $b(x)=1$),
$\mathbb{P} \{\aleph(\mathbf{p}_{50})=i \}$ for
$i=1,\ldots,50$. %
\textit{First Column}: Partition sizes; %
\textit{Second Column}: True probability of partition sizes; %
\textit{Third Column}: Mean of 10,000 approximates of the probabilities
according to algorithm A.1; %
\textit{Fourth Column}: Standard error of 10,000 approximates of the
probabilities according to algorithm A.1; %
\textit{Fifth Column}: Mean of 10,000 approximates of the probabilities
according to algorithm A.2; %
\textit{Sixth Column}: Standard error of 10,000 approximates of the
probabilities according to algorithm A.2; %
\textit{Seventh Column}: Mean of 10,000 approximates of the probabilities
according to algorithm A.3; %
\textit{Eighth Column}: Standard error of 10,000 approximates of the
probabilities according to algorithm A.3; %
\textit{Ninth Column}: Mean of 10,000 approximates of the probabilities
according to algorithm A.4; %
\textit{Tenth Column}: Standard error of 10,000 approximates of the
probabilities according to algorithm A.4}\label{tab2}
\begin{tabular*}{\textwidth}{@{\extracolsep{\fill}}llllllllll@{}}
\hline
& $\mathbb{P} \{\aleph(\mathbf{p}_n)=i \}$ & \multicolumn{2}{l}{{A.1 algorithm}}
& \multicolumn{2}{l}{ {A.2 algorithm}} & \multicolumn{2}{l}{{A.3 algorithm}} &
\multicolumn{2}{l}{{A.4 algorithm}} \\[-5pt]
&& \multicolumn{2}{l}{\hrulefill} & \multicolumn{2}{l}{\hrulefill} &
\multicolumn{2}{l}{\hrulefill} & \multicolumn{2}{l}{\hrulefill}\\
& {True} &  {Mean of the} &  {SE of the} &  {Mean of the} &
{SE of the} &  {Mean of the} &  {SE of the} &  {Mean of the} &  {SE of the} \\
$i$ &  {probability} &  {approximates} &  {approximates} & {approximates} &  {approximates}
&  {approximates} &  {approximates} &  {approximates} &  {approximates} \\
\hline
\hphantom{0}1 & 0.000035 & 0.000036 & 0.000059 & 0.000035 & 0.000060 & 0.000035 &
0.000061 & 0.000035 & 0.000061 \\
\hphantom{0}2 & 0.000291 & 0.000294 & 0.000170 & 0.000290 & 0.000170 & 0.000288 &
0.000178 & 0.000293 & 0.000183 \\
\hphantom{0}3 & 0.001234 & 0.001227 & 0.000353 & 0.001238 & 0.000356 & 0.001229 &
0.000389 & 0.001241 & 0.000396 \\
\hphantom{0}4 & 0.003621 & 0.003621 & 0.000599 & 0.003624 & 0.000604 & 0.003620 &
0.000701 & 0.003622 & 0.000703 \\
\hphantom{0}5 & 0.008265 & 0.008274 & 0.000916 & 0.008266 & 0.000909 & 0.008249 &
0.001101 & 0.008262 & 0.001083 \\
\hphantom{0}6 & 0.015684 & 0.015709 & 0.001237 & 0.015683 & 0.001250 & 0.015647 &
0.001554 & 0.015666 & 0.001551 \\
\hphantom{0}7 & 0.025793 & 0.025807 & 0.001598 & 0.025753 & 0.001594 & 0.025786 &
0.001989 & 0.025825 & 0.001997 \\
\hphantom{0}8 & 0.037832 & 0.037827 & 0.001887 & 0.037807 & 0.001896 & 0.037806 &
0.002365 & 0.037843 & 0.002367 \\
\hphantom{0}9 & 0.050538 & 0.050540 & 0.002207 & 0.050498 & 0.002233 & 0.050500 &
0.002645 & 0.050583 & 0.002657 \\
10 & 0.062458 & 0.062450 & 0.002393 & 0.062502 & 0.002396 & 0.062445 &
0.002807 & 0.062497 & 0.002823 \\
11 & 0.072282 & 0.072260 & 0.002591 & 0.072270 & 0.002557 & 0.072238 &
0.002921 & 0.072246 & 0.002880 \\
12 & 0.079070 & 0.079074 & 0.002674 & 0.079123 & 0.002698 & 0.079058 &
0.002932 & 0.079112 & 0.002914 \\
13 & 0.082360 & 0.082371 & 0.002745 & 0.082353 & 0.002724 & 0.082412 &
0.002945 & 0.082404 & 0.002924 \\
14 & 0.082159 & 0.082126 & 0.002758 & 0.082176 & 0.002764 & 0.082174 &
0.002894 & 0.082143 & 0.002915 \\
15 & 0.078842 & 0.078821 & 0.002708 & 0.078836 & 0.002675 & 0.078878 &
0.002847 & 0.078786 & 0.002846 \\
16 & 0.073037 & 0.072990 & 0.002618 & 0.073001 & 0.002614 & 0.073106 &
0.002783 & 0.073063 & 0.002804 \\
17 & 0.065486 & 0.065510 & 0.002469 & 0.065512 & 0.002465 & 0.065489 &
0.002706 & 0.065455 & 0.002707 \\
18 & 0.056942 & 0.056976 & 0.002343 & 0.056984 & 0.002316 & 0.056960 &
0.002585 & 0.056956 & 0.002560 \\
19 & 0.048085 & 0.048085 & 0.002165 & 0.048092 & 0.002160 & 0.048113 &
0.002429 & 0.048047 & 0.002398 \\
20 & 0.039472 & 0.039478 & 0.001975 & 0.039484 & 0.001980 & 0.039473 &
0.002225 & 0.039428 & 0.002223 \\
21 & 0.031516 & 0.031532 & 0.001743 & 0.031496 & 0.001741 & 0.031529 &
0.002037 & 0.031498 & 0.002037 \\
\hline
\end{tabular*}
\end{sidewaystable*}
\setcounter{table}{1}

\begin{sidewaystable*}
\tablewidth=\textwidth
\tabcolsep=0pt
\caption{(Continued)}
\begin{tabular*}{\textwidth}{@{\extracolsep{\fill}}llllllllll@{}}
\hline
& $\mathbb{P} \{\aleph(\mathbf{p}_n)=i \}$ & \multicolumn{2}{l}{{A.1 algorithm}}
& \multicolumn{2}{l}{ {A.2 algorithm}} & \multicolumn{2}{l}{{A.3 algorithm}} &
\multicolumn{2}{l}{{A.4 algorithm}} \\[-5pt]
&& \multicolumn{2}{l}{\hrulefill} & \multicolumn{2}{l}{\hrulefill} &
\multicolumn{2}{l}{\hrulefill} & \multicolumn{2}{l}{\hrulefill}\\
& {True} &  {Mean of the} &  {SE of the} &  {Mean of the} &
{SE of the} &  {Mean of the} &  {SE of the} &  {Mean of the} &  {SE of the} \\
$i$ &  {probability} &  {approximates} &  {approximates} & {approximates} &  {approximates}
&  {approximates} &  {approximates} &  {approximates} &  {approximates} \\
\hline
22 & 0.024480 & 0.024454 & 0.001541 & 0.024485 & 0.001554 & 0.024485 &
0.001765 & 0.024466 & 0.001769 \\
23 & 0.018499 & 0.018503 & 0.001338 & 0.018488 & 0.001361 & 0.018484 &
0.001535 & 0.018502 & 0.001544 \\
24 & 0.013595 & 0.013589 & 0.001155 & 0.013583 & 0.001166 & 0.013584 &
0.001308 & 0.013604 & 0.001307 \\
25 & 0.009711 & 0.009716 & 0.000978 & 0.009696 & 0.000982 & 0.009703 &
0.001091 & 0.009714 & 0.001110 \\
26 & 0.006738 & 0.006750 & 0.000824 & 0.006743 & 0.000819 & 0.006738 &
0.000907 & 0.006725 & 0.000910 \\
27 & 0.004537 & 0.004544 & 0.000670 & 0.004547 & 0.000676 & 0.004544 &
0.000738 & 0.004538 & 0.000732 \\
28 & 0.002960 & 0.002961 & 0.000545 & 0.002959 & 0.000547 & 0.002958 &
0.000593 & 0.002960 & 0.000585 \\
29 & 0.001870 & 0.001866 & 0.000434 & 0.001867 & 0.000432 & 0.001870 &
0.000460 & 0.001872 & 0.000463 \\
30 & 0.001141 & 0.001144 & 0.000337 & 0.001144 & 0.000338 & 0.001136 &
0.000354 & 0.001146 & 0.000359 \\
31 & 0.000672 & 0.000672 & 0.000260 & 0.000674 & 0.000257 & 0.000667 &
0.000269 & 0.000673 & 0.000273 \\
32 & 0.000381 & 0.000379 & 0.000195 & 0.000380 & 0.000195 & 0.000381 &
0.000203 & 0.000384 & 0.000203 \\
33 & 0.000207 & 0.000207 & 0.000144 & 0.000206 & 0.000143 & 0.000208 &
0.000148 & 0.000205 & 0.000147 \\
34 & 0.000108 & 0.000108 & 0.000105 & 0.000108 & 0.000104 & 0.000107 &
0.000105 & 0.000107 & 0.000105 \\
35 & 0.000054 & 0.000053 & 0.000073 & 0.000053 & 0.000072 & 0.000054 &
0.000075 & 0.000053 & 0.000073 \\
36 & 0.000025 & 0.000025 & 0.000050 & 0.000025 & 0.000051 & 0.000026 &
0.000052 & 0.000025 & 0.000050 \\
37 & 0.000011 & 0.000011 & 0.000033 & 0.000012 & 0.000034 & 0.000011 &
0.000033 & 0.000012 & 0.000034 \\
38 & 0.000005 & 0.000005 & 0.000023 & 0.000005 & 0.000022 & 0.000005 &
0.000022 & 0.000005 & 0.000022 \\
39 & 0.000002 & 0.000002 & 0.000014 & 0.000002 & 0.000014 & 0.000002 &
0.000013 & 0.000002 & 0.000014 \\
40 & 0.000001 & 0.000001 & 0.000008 & 0.000001 & 0.000009 & 0.000001 &
0.000009 & 0.000001 & 0.000008 \\
41 & 0.000000 & 0.000000 & 0.000006 & 0.000000 & 0.000005 & 0.000000 &
0.000005 & 0.000000 & 0.000005 \\
42 & 0.000000 & 0.000000 & 0.000003 & 0.000000 & 0.000002 & 0.000000 &
0.000002 & 0.000000 & 0.000002 \\
43 & 0.000000 & 0.000000 & 0.000001 & 0.000000 & 0.000001 & 0.000000 &
0.000001 & 0.000000 & 0.000001 \\
44 & 0.000000 & 0.000000 & 0.000000 & 0.000000 & 0.000000 & 0.000000 &
0.000000 & 0.000000 & 0.000000 \\
45 & 0.000000 & 0.000000 & 0.000000 & 0.000000 & 0.000000 & 0.000000 &
0.000000 & 0.000000 & 0.000000 \\
46 & 0.000000 & 0.000000 & 0.000000 & 0.000000 & 0.000000 & 0.000000 &
0.000000 & 0.000000 & 0.000000 \\
47 & 0.000000 & 0.000000 & 0.000000 & 0.000000 & 0.000000 & 0.000000 &
0.000000 & 0.000000 & 0.000000 \\
48 & 0.000000 & 0.000000 & 0.000000 & 0.000000 & 0.000000 & 0.000000 &
0.000000 & 0.000000 & 0.000000 \\
49 & 0.000000 & 0.000000 & 0.000000 & 0.000000 & 0.000000 & 0.000000 &
0.000000 & 0.000000 & 0.000000 \\
50 & 0.000000 & 0.000000 & 0.000000 & 0.000000 & 0.000000 & 0.000000 &
0.000000 & 0.000000 & 0.000000 \\
\hline
\end{tabular*}
\end{sidewaystable*}

Finally, we note that the applications of Bayesian non-parametric mixture
models in Bayesian statistics is steadily increasing (see
Lo \cite{Lo1984}; James \cite{James2002}; Ishwaran and James \cite{IJ2001,IJ2003};
James \textit{et al.} \cite{JLP2008}; Lijoi \textit{et al.} \cite{LMP2005a,LMP2005b,LMP2007b}). In
particular, time series model mixing over random probability
measures has been considered recently in Griffin and Steel \cite{GS2006},
Lau and So \cite{LS2008,LS2011} and Lau and Cripps \cite{LC2012}. Often in mixture models
over the normalized tilted \textsc{crm}, it is necessary to consider a
collection of latent variables, which is a sample of the normalized
tilted \textsc{crm} but these latent variables are not observed
directly. In fact, sampling latent variables is essential for
approximating estimates of parameters of interest that are functions
of the latent variables, $\mathbf{X}_n$ because of the high
cardinality of the posterior distribution due to combinatorial
property of the latent variables. As a result, sampling schemes for
$\mathbf{X}_n$ are required for estimation based on the conditional
Blackwell--MacQueen P\'olya urn formula and the distributions of
$\mathbf{X}_n$ and partitions $\mathbf{p}_n$. In this article, we
have provided the marginal distributions of $\mathbf{X}_n$ and
$\mathbf{p}_n$, both conditional and unconditional on $U_n$, that
are essential elements in implementing mixture models over the
normalized tilted \textsc{crm}.

\begin{appendix}

\section*{Appendix}\label{app}
\subsection{\texorpdfstring{Proof of Theorem \protect\ref{thm1posteriorG}}
{Proof of Theorem 3.1}}\label{thm1proof}

\setcounter{equation}{0}

Consider the joint distribution of $(\widetilde\mu_t,\mathbf{X}_n)$,
which is given by (\ref{jointlaw}) and (\ref{expression1}), that is
%
\begin{equation}
\label{Thm1eqn1} \mathcal{P}_{\widetilde{\mu}_t,\mathbf{X}_n}(A,B_1,
\ldots,B_n)%
=
\frac{1}%
{E[h^\prime(\widetilde{
\mu}(\mathbb{X}))\widetilde{\mu}(\mathbb{X})^{-q}]} \int
_{A}\frac{h^\prime(\mu(\mathbb{X}))}{\mu(\mathbb{X})^{n+q}}%
\Biggl(\prod
_{i=1}^n \mu(B_i) \Biggr)
\mathcal{P}_{\widetilde{\mu}}(\mathrm{d}\mu) .%
\end{equation}
In (\ref{Thm1eqn1}), the \textsc{crm} $\widetilde\mu(\cdot)$ is
now replaced by the linear functional of the Poisson process
$\int_{\mathbb{R}^{+}}sN(\cdot,\mathrm{d}s)$ on $\mathbb{X}$ (\ref{CRM}).
Following by writing the right hand side of (\ref{Thm1eqn1}) without
the integrals, that is given by
%
\begin{eqnarray}
\label{Thm1eqn2} 
\frac{h^\prime(\int_{\mathbb{X}\times\mathbb{R}^{+}} s N(\mathrm{d}x,\mathrm{d}s))}{(\int_{\mathbb{X}\times\mathbb{R}^{+}} s N(\mathrm{d}x,\mathrm{d}s))^{n+q}}%
\Biggl(\prod
_{i=1}^n 
s_i
N(\mathrm{d}x_i,\mathrm{d}s_i) \Biggr) \mathcal{P}_{N}(\mathrm{d}N)
,%
\end{eqnarray}
where this Poisson process $N$ has distribution denoted by
$\mathcal{P}_{N}$, has the intensity measure $\nu$ same as that of
$\widetilde\mu$ (\ref{CRM}), and belong to the set in
$\mathcal{B}(\mathcal{M}_{\mathbb{X}\times\mathbb{R}^+})$
corresponding to the set $A$ that $\mu\in
A\in\mathcal{B}(\mathcal{M}_{\mathbb{X}})$.
Here $(x_i,s_i)$ for $i=1,\ldots,n$ represent the points generated
from the Poisson process $N$. For $k=1,\ldots,\aleph(\mathbf{p}_n)$,
$(Y_k,J_k)$ denotes the distinct points of the $(x_i,s_i)$'s and
$\aleph(\mathbf{p}_n)$ denotes the number of the distinct points. In
addition, we take $ \{J_1,\ldots,J_{\aleph(\mathbf{p}_n)} \}$
to be the augmented variables. 
We then apply the Fubini theorem following from an application of
Lemma 2.2 of James \cite{James2002}, page 8 (see also James \cite{James2002}
for some detail discussions), to yield the joint distribution of $N$,
$\mathbf{p}_n$ and $(Y_k,J_k)$ for
$k=1,\ldots,\aleph(\mathbf{p}_n)$, which is given by
(\ref{Thm1eqn2}) without the proportional constant,
%
\begin{eqnarray}
\label{Thm1eqn3} 
\frac{h^\prime(\int_{\mathbb{X}\times\mathbb{R}^{+}} s N^{(\mathbf
{X}_n)}(\mathrm{d}x,\mathrm{d}s))}{(\int_{\mathbb{X}\times\mathbb{R}^{+}} s N^{(\mathbf
{X}_n)}(\mathrm{d}x,\mathrm{d}s))^{n+q}}%
\mathcal{P}_{N}(\mathrm{d}N)%
\prod_{k=1}^{\aleph(\mathbf{p}_n)}
\eta(\mathrm{d}Y_k) s_k^{n_k}\rho_{Y_k}(\mathrm{d}s_k)
,
\end{eqnarray}
where
$N^{(\mathbf{X}_n)}\overset{d}{=}N+\sum_{k=1}^{\aleph(\mathbf
{p}_n)}\delta_{(Y_k,J_k^{(\mathbf{X}_n)})}$.
Here the distribution of $N^{(\mathbf{X}_n)}$ is identical to the
distribution of the sum over a Poisson process $N$ and the fixed
points of discontinuity at $(Y_k,J_k^{(\mathbf{X}_n)})$. The Poisson
process $N$ has the intensity measure $\nu$ and be independent of
$(Y_k,J_k^{(\mathbf{X}_n)})$
for $k=1,\ldots,\aleph(\mathbf{p}_n)$. The pairs
$(Y_k,J_k^{(\mathbf{X}_n)})$ for
$k=1,\ldots,\aleph(\mathbf{p}_n)$ and $\mathbf{p}_n$ has the joint
distribution $\prod_{k=1}^{\aleph(\mathbf{p}_n)}
[\eta(\mathrm{d}Y_k)\times\mathbb{P}(J_k^{(\mathbf{X}_n)}\in
\mathrm{d}s_k)\times\tau_{n_k,Y_k}(0) ]$ where
$\mathbb{P}(J_k^{(\mathbf{X}_n)}\in \mathrm{d}s_k)=s_k^{n_k}\rho_{Y_k}(\mathrm{d}s_k)
/ \int_{\mathbb{R}^+}s_k^{n_k}\rho_{Y_k}(\mathrm{d}s_k)$ and
$\tau_{n_k,Y_k}(0)$ defined in (\ref{tau}).
The distribution of
$N^{(\mathbf{X}_n)}$ now is denoted by
$\mathcal{P}_{N}^{(\mathbf{X}_n)}$, so (\ref{Thm1eqn3}) could be
reduced to
%
\begin{eqnarray}
\label{Thm1eqn4} \frac{h^\prime(\int_{\mathbb{X}\times\mathbb{R}^{+}}sN(\mathrm{d}x,\mathrm{d}s))}%
{(\int_{\mathbb{X}\times\mathbb{R}^{+}}sN(\mathrm{d}x,\mathrm{d}s)
)^{n+q}}%
\mathcal{P}_{N}^{(\mathbf{X}_n)}(\mathrm{d}N)
.
\end{eqnarray}
The Poisson linear functional appears in (\ref{Thm1eqn4}),
$\int_{\mathbb{R}^{+}}sN^{(\mathbf{X}_n)}(\mathrm{d}x,\mathrm{d}s)$
($=:\widetilde\mu^{(\mathbf{X}_n)}(\mathrm{d}x)$), is a completely random
measure according to (\ref{representation}) (see also
Daley and Vere-Jones \cite{DV2008}, Theorem 10.1.III, page 79), such that
%
\begin{eqnarray}
\label{PosteriorCPM} \widetilde\mu^{(\mathbf{X}_n)}(\mathrm{d}x) \overset{d} {=}\int
_{\mathbb{R}^{+}}s N (\mathrm{d}x,\mathrm{d}s) +\sum_{k=1}^{\aleph(\mathbf{p}_n)}
J_k^{(\mathbf{X}_n)} \delta_{Y_k}(\mathrm{d}x) \overset{d} {=}
\widetilde\mu(\mathrm{d}x) +\sum_{k=1}^{\aleph(\mathbf{p}_n)}
J_k^{(\mathbf{X}_n)}\delta_{Y_k}(\mathrm{d}x) ,
\end{eqnarray}
where $\widetilde\mu$ denote a \textsc{crm} (\ref{CRM}) derived from
the Poisson process $N$. So, the total mass
$\widetilde\mu^{(\mathbf{X}_n)}(\mathbb{X})$ in (\ref{Thm1eqn4}) is
given by
%
\begin{eqnarray}
\label{Thm1eqn5} \widetilde\mu^{(\mathbf{X}_n)}(\mathbb{X})\overset{d} {=}
\widetilde\mu (\mathbb{X}) +\sum_{k=1}^{\aleph(\mathbf{p}_n)}
J_k^{(\mathbf{X}_n)} .
\end{eqnarray}
This immediately reveals the distribution of the tilted completely
random measure, $\widetilde\mu^{(\mathbf{X}_n)}_t$ (Definition
\ref{Deftiltedmeasure}), that is
$\mathcal{P}_{\widetilde\mu_t}^{(\mathbf{X}_n)}(\mathrm{d}\mu)\propto
\frac{h^\prime(\mu^{(\mathbf{X}_n)}(\mathbb{X}))}%
{(\mu^{(\mathbf{X}_n)}(\mathbb{X}))^{n+q}}%
\mathcal{P}_{\widetilde\mu}^{(\mathbf{X}_n)}(\mathrm{d}\mu)$.
Lastly, the distributional identity derived from the a sample
$\widetilde\mu^{(\mathbf{X}_n)}_t$ from $\mathcal{P}_{\widetilde\mu
_t}^{(\mathbf{X}_n)}$,
$G^{(\mathbf{X}_n)}(\cdot) \overset{d}=
\widetilde\mu^{(\mathbf{X}_n)}_t(\cdot)/\widetilde\mu^{(\mathbf
{X}_n)}_t(\mathbb{X})$
on $\mathbb{X}$ is obtained by the same arguments in the proof of
Theorem 2 of James \textit{et al.} \cite{JLP2008}, page 96. Thus, the proof is complete.

\subsection{\texorpdfstring{Proof of Theorem \protect\ref{thm2conditionalposteriorG}}
{Proof of Theorem 3.2}}\label{thm2proof}

Following the proof of Theorem \ref{thm1posteriorG} from the
beginning to (\ref{Thm1eqn3}) with chosen $h^\prime(x)=\mathrm{e}^{-\gamma x}$.
Using the facts from (\ref{PosteriorCPM}) and (\ref{Thm1eqn5}), the
distribution (\ref{Thm1eqn3}) becomes
%
\begin{eqnarray}
\label{Thm2eqn1} \frac{\mathrm{e}^{-\gamma(\mu(\mathbb{X})+\sum_{k=1}^{\aleph(\mathbf{p}_n)} s_k
)}}%
{(\mu(\mathbb{X})+\sum
_{k=1}^{\aleph(\mathbf{p}_n)} s_k
)^{n+q}}%
\mathcal{P}_{\widetilde\mu}(\mathrm{d}\mu) \prod
_{k=1}^{\aleph(\mathbf{p}_n)} \eta(\mathrm{d}Y_k)
s_k^{n_k} \rho_{Y_k}(\mathrm{d}s_k )
.%
\end{eqnarray}
This is the joint distribution of $\widetilde\mu$, $\mathbf{p}_n$ and
$Y_k$ for $k=1,\ldots,\aleph(\mathbf{p}_n)$.
Then we apply the gamma identity (\ref{gammaidentity}) on the term
$(\mu(\mathbb{X})+\sum_{k=1}^{\aleph(\mathbf{p}_n)} s_k)^{-(n+q)}$
in (\ref{Thm2eqn1})
%
\begin{eqnarray}
\label{Thm2eqn2} \frac{1}{(\mu(\mathbb{X})+\sum_{k=1}^{\aleph(\mathbf{p}_n)} s_k)^{n+q}}%
=\frac{1}{\Gamma(n+q)}\int
_{\mathbb{R}^{+}}%
\mathrm{e}^{-(\mu(\mathbb{X})+\sum_{k=1}^{\aleph(\mathbf{p}_n)}
s_k)u}u^{n+q-1}\,\mathrm{d}u
,%
\end{eqnarray}
Here the augmented variable $U_n$ is introduced due to the gamma
identity. Incorporating (\ref{Thm2eqn1}) with (\ref{Thm2eqn2}) and
omiting the integrals, it turns out that the following is the joint
distribution of $U_n$, $\widetilde\mu$, $\mathbf{p}_n$ and $Y_k$ for
$k=1,\ldots,\aleph(\mathbf{p}_n)$,
%
\begin{eqnarray}
\label{Thm2eqn3} \mathrm{e}^{-(\gamma+u)\mu(\mathbb{X})} \mathcal{P}_{\widetilde\mu}(\mathrm{d}\mu)
\frac{u^{n+q-1}}{\Gamma(n+q)} \prod_{k=1}^{\aleph(\mathbf{p}_n)}
\eta(\mathrm{d}Y_k) s_k^{n_k} \mathrm{e}^{-s_k(\gamma+u)}%
\rho_{Y_k}(\mathrm{d}s_k)%
\,\mathrm{d}u .%
\end{eqnarray}
The disintegration between terms in $\mathrm{e}^{-(\gamma+u)\mu(\mathbb{X})}$
and $\mathcal{P}_{\widetilde\mu}(\mathrm{d}\mu)$ in (\ref{Thm2eqn3}) yields
$\mathrm{e}^{-\psi_0(\gamma+u)}$ and
$\mathcal{P}_{\widetilde\mu^{(u)}}(\mathrm{d}\mu)$ where $\psi_0$ is defined
in (\ref{psi}) and the completely random measure
$\widetilde\mu^{(u)}$ has the intensity measure
$\nu^{(u)}(\mathrm{d}x,\mathrm{d}s)=\eta(\mathrm{d}x)\mathrm{e}^{-(\gamma+u)s}\rho_x(\mathrm{d}s)$. Then,
(\ref{Thm2eqn3}) turns out to be
%
\begin{eqnarray}
\label{Thm2eqn4} \mathcal{P}_{\widetilde\mu^{(u)}}(\mathrm{d}\mu)%
\mathrm{e}^{-\psi_0(\gamma+u)}
\frac{u^{n+q-1}}{\Gamma(n+q)} \prod_{k=1}^{\aleph(\mathbf{p}_n)}
\eta(\mathrm{d}Y_k) s_k^{n_k} \mathrm{e}^{-s_k(\gamma+u)}%
\rho_{Y_k}(\mathrm{d}s_k)%
\, \mathrm{d}u .%
\end{eqnarray}
So, conditional on $U_n$ and $\mathbf{X}_n$, the process
$\widetilde\mu^{(U_n)}+\sum_{k=1}^{\aleph(\mathbf
{p}_n)}J_k^{(U_n,\mathbf{X}_n)}\delta_{Y_k}$
is a completely random measure (\ref{representation}) (see also
Daley and Vere-Jones \cite{DV2008}, page 79, Theorem 10.1.III) such that
%
\begin{eqnarray}
\widetilde\mu^{(U_n,\mathbf{X}_n)}(\mathrm{d}x)%
\overset{d} {=}%
\widetilde\mu^{(U_n)}(\mathrm{d}x)+%
\sum_{k=1}^{\aleph(\mathbf{p}_n)}
J_k^{(U_n,\mathbf{X}_n)}\delta_{Y_k}(\mathrm{d}x) ,
\end{eqnarray}
where each $J_k^{(U_n,\mathbf{X}_n)}$ has the density
$\mathbb{P}(J_k^{(U_n,\mathbf{X}_n)}\in \mathrm{d}s) = \frac{s^{n_k}
\mathrm{e}^{-s(\gamma+U_n)}\rho_{Y_k}(\mathrm{d}s)}{\int_{\mathbb{R}^+} s^{n_k}
\mathrm{e}^{-s(\gamma+U_n)}\rho_{Y_k}(\mathrm{d}s)}$ for
$k=1,\ldots,\allowbreak \aleph(\mathbf{p}_n)$ and
$ \{J_1^{(U_n,\mathbf{X}_n)},\ldots,J_{\aleph(\mathbf
{p}_n)}^{(U_n,\mathbf{X}_n)} \}$
are conditionally independent. Furthermore, conditional on $U_n$,
$\widetilde\mu^{(U_n)}$ and
$ \{J_1^{(U_n,\mathbf{X}_n)},\ldots,J_{\aleph(\mathbf
{p}_n)}^{(U_n,\mathbf{X}_n)} \}$
are independent. Lastly the distributional identity
$\widetilde{G}^{(U_n,\mathbf{X}_n)}(\cdot)\overset{d}{=}\widetilde\mu^{(U_n,
\mathbf{X}_n)}(\cdot)/\widetilde\mu^{(U_n,\mathbf{X}_n)}(\mathbb{X})$ is obtained by the
same arguments in the proof of Theorem 2 of James \textit{et al.} \cite{JLP2008}, page
96, and $\widetilde G^{(U_n,\mathbf{X}_n)}$ is a normalized completely
random measure. Thus, the proof is complete.

\subsection{\texorpdfstring{Proof of Proposition \protect\ref{prop1BMUrnpredictivedistribution}}
{Proof of Proposition 3.1}}\label{prop1proof} Following the definition of the
predictive distribution,
%
\begin{eqnarray}
\label{Prop1eqn1} \mathbb{P} \{X_{n+1}\in \mathrm{d}x|\mathbf{X}_{n}
\}%
=\int_{\mathbb{R}^{+}} E \bigl[\widetilde
G^{(U_n,\mathbf{X}_n)}(\mathrm{d}x)\vert \mathbf{X}_{n},U_{n}=u \bigr]
f_{U_{n}\vert\mathbf{X}_{n}}(u)\,\mathrm{d}u ,%
\end{eqnarray}
and using the result from Theorem \ref{thm2conditionalposteriorG}, the expectation (\ref{Prop1eqn1}) inside the integral is given
by
%
\begin{eqnarray}
\label{Prop1eqn2} && E \bigl[\widetilde G^{(U_n,\mathbf{X}_n)}(\mathrm{d}x)\vert
\mathbf{X}_{n},U_{n} \bigr]\nonumber\\
&&\quad= E \biggl[\frac{\widetilde\mu^{(U_n)}(\mathrm{d}x)}{\widetilde\mu^{(U_n)}(\mathbb
{X})+\sum_{k^{\prime}=1}^{\aleph(\mathbf{p}_{n})}J_{k^{\prime
}}^{(U_n,\mathbf{X}_n)}}
\Big\vert\mathbf{X}_{n},U_{n} \biggr]%
\\
&&\qquad{}+ %
\sum_{k=1}^{\aleph(\mathbf{p}_{n})} E \biggl[
\frac{J_{k}^{(U_n,\mathbf{X}_n)}}{\widetilde\mu^{(U_n)}(\mathbb
{X})+\sum_{k^{\prime}=1}^{\aleph(\mathbf{p}_{n})}J_{k^{\prime
}}^{(U_n,\mathbf{X}_n)}} \Big\vert\mathbf{X}_{n},U_{n}
\biggr]%
\delta_{Y_k}(\mathrm{d}x) .
\nonumber
\end{eqnarray}
We use Theorem \ref{thm2conditionalposteriorG} to obtain the
explicit results of the expectations (\ref{Prop1eqn2}) according to
the conditional distributions of $\widetilde\mu^{(U_n)}$ and
$J_{k}^{(U_n,\mathbf{X}_n)}$ for $k=1,\ldots,\aleph(\mathbf{p}_n)$.
There are two expected values (\ref{Prop1eqn2}) considered in the
following. Firstly, an exponential identity in the first integral of
(\ref{Prop1eqn2}) yields,
%
\begin{eqnarray}
\label{Prop1eqn3} 
&& E \biggl[\frac{\widetilde\mu^{(U_n)}(\mathrm{d}x)}{\mu^{(U_n)}(\mathbb{X})+\sum_{k^{\prime}=1}^{\aleph(\mathbf{p}_{n})}J_{k^{\prime}}^{(U_n,\mathbf
{X}_n)}} \Big\vert
\mathbf{X}_{n},U_{n} \biggr]\nonumber\\[-8pt]\\[-8pt]
&&\quad
=\int
_{\mathbb{R}^{+}}E \Biggl[\widetilde\mu^{(U_n)}(\mathrm{d}x)\mathrm{e}^{-z\widetilde\mu
^{(U_n)}(\mathbb{X})}
\prod_{k^{\prime}=1}^{\aleph(\mathbf{p}_{n})}\mathrm{e}^{-zJ_{k^{\prime
}}^{(U_n,\mathbf{X}_n)}} \Big\vert
\mathbf{X}_{n},U_{n} \Biggr]\mathrm{d}z. \nonumber%
\end{eqnarray}
%
We apply a change of measure on $\widetilde\mu^{(U_n)}$ and compute
directly on the expectation with respect to
$J_{k}^{(U_n,\mathbf{X}_n)}$ for $k=1,\ldots,\aleph(\mathbf{p}_n)$,
(\ref{Prop1eqn3}) turns out to be
%
\begin{eqnarray}
\label{Prop1eqn5} \eta(\mathrm{d}x)\times\int_{\mathbb{R}^{+}}\mathrm{e}^{-\psi_{\gamma+U_{n}}(z)}
\tau_{1,x}(\gamma+z+U_{n})%
\prod
_{k=1}^{\aleph(\mathbf{p}_{n})}%
\frac{\tau_{n_{k},Y_k}(\gamma+z+U_{n})}%
{
\tau_{n_{k},Y_k}(\gamma+U_{n})}\,\mathrm{d}z .%
\end{eqnarray}
Marginalizing over $U_n$ on (\ref{Prop1eqn5}) is given by
%
\begin{eqnarray}
\label{Prop1eqn6} &&\int_{\mathbb{R}^{+}} E \biggl[\frac{\widetilde\mu^{(U_n)}(\mathrm{d}x)}{\widetilde\mu^{(U_n)}(\mathbb
{X})+\sum_{k^{\prime}=1}^{\aleph(\mathbf{p}_{n})}J_{k^{\prime
}}^{(U_n,\mathbf{X}_n)}}
\vert\mathbf{X}_{n},U_n \biggr]%
f_{U_n\vert\mathbf{X}_n}(u)\,\mathrm{d}u
\nonumber\\
&&\quad=\eta(\mathrm{d}x)\times\int_{\mathbb{R}^{+}\times\mathbb{R}^{+}} \mathrm{e}^{-\psi_{\gamma+u}(z)}%
\tau_{1,x}(\gamma+z+u)
\\
&&\hspace*{86pt}{}\times\frac{%
\mathrm{e}^{-\psi_{0}(\gamma+u)}%
\prod_{k=1}^{\aleph(\mathbf{p}_{n})}%
\tau_{n_{k},Y_{k}}(\gamma+z+u)u^{n+q-1}}%
{\int_{\mathbb{R}^{+}}\mathrm{e}^{-\psi_{0}(\gamma+u)}%
\prod_{k=1}^{\aleph(\mathbf{p}_{n})}\tau_{n_{k},Y_{k}}(\gamma
+u)u^{n+q-1}\,\mathrm{d}u}\,\mathrm{d}u\,\mathrm{d}z.
\nonumber
\end{eqnarray}
Taking the equality
$\psi_{0}(\gamma+z+u)=\psi_{\gamma+u}(z)+\psi_{0}(\gamma+u)$ and
transforming the upper integral with $(w,v)=(u+z,u)$,
(\ref{Prop1eqn6}) becomes
\begin{eqnarray*}
\label{Prop1eqn7} &&\eta(\mathrm{d}x)\times\frac{\int_{\mathbb{R}^{+}}%
\mathrm{e}^{-\psi_{0}(\gamma+w)}%
\tau_{1,x}(\gamma+w)\prod_{k=1}^{\aleph(\mathbf{p}_{n})}%
\tau_{n_{k},Y_{k}}(\gamma+w) {w^{n+q}}/({n+q})\,\mathrm{d}w}%
{\int
_{\mathbb{R}^{+}}\mathrm{e}^{-\psi_{0}(\gamma+u)}%
\prod
_{k=1}^{\aleph(\mathbf{p}_{n})}\tau_{n_{k},Y_{k}}(\gamma
+u)u^{n+q-1}\,\mathrm{d}u}%
\\
&&\quad=\eta(\mathrm{d}x)\times\frac{1}{n+q}\int_{\mathbb{R}^{+}}%
\tau_{1,x}(\gamma+u)uf_{U_n\vert\mathbf{X}_n}(u)\,\mathrm{d}u .
\end{eqnarray*}
Now consider the second expectation of (\ref{Prop1eqn2}). As before,
an exponential identity in the second integral yields,
%
\begin{eqnarray}
\label{Prop1eqn8} &&
E \biggl[\frac{J_{k}^{(U_n,\mathbf{X}_n)}}{\mu^{(U_n)}(\mathbb{X})+\sum_{k^{\prime}=1}^{\aleph(\mathbf{p}_{n})}J_{k^{\prime}}^{(U_n,\mathbf
{X}_n)}} \Big\vert
\mathbf{X}_{n},U_{n} \biggr]%
\nonumber\\[-8pt]\\[-8pt]
&&\quad
=\int
_{\mathbb{R}^{+}}E \Biggl[J_{k}^{(U_n,\mathbf{X}_n)}\mathrm{e}^{-z\mu
^{(U_n)}(\mathbb{X})}
\prod_{k^{\prime}=1}^{\aleph(\mathbf{p}_{n})}\mathrm{e}^{-zJ_{k^{\prime
}}^{(U_n,\mathbf{X}_n)}} \Big\vert
\mathbf{X}_{n},U_{n} \Biggr]\,\mathrm{d}z . \nonumber%
\end{eqnarray}
By direct computation, (\ref{Prop1eqn8}) turns out to be
%
\begin{eqnarray}
\label{Prop1eqn10} \int_{\mathbb{R}^{+}}%
\mathrm{e}^{-\psi_{\gamma+U_{n}}(z)}%
\frac{\tau_{n_{k}+1,Y_{k}}(\gamma+z+U_{n})}{\tau_{n_{k},Y_{k}}(\gamma
+U_{n})}%
\prod_{k^{\prime}=1,k^{\prime}\neq k}^{\aleph(\mathbf{p}_{n})}%
\frac{\tau_{n_{k^{\prime}},Y_{k^{\prime}}}(\gamma+z+U_{n})}%
{\tau_{n_{k^{\prime}},Y_{k^{\prime}}}(\gamma+U_{n})}%
\,\mathrm{d}z .%
\end{eqnarray}
Marginalizing over $U_n$ on (\ref{Prop1eqn10}) is given by
\begin{eqnarray}
\label{Prop1eqn11} && \int_{\mathbb{R}^{+}}%
E \biggl[
\frac{J_{k}^{(U_n,\mathbf{X}_n)}}{\widetilde\mu^{(U_n)}(\mathbb
{X})+\sum_{k^{\prime}=1}^{\aleph(\mathbf{p}_{n})}J_{k^{\prime
}}^{(U_n,\mathbf{X}_n)}}\Big \vert\mathbf{X}_{n},U_n
\biggr]%
f_{U_n\vert\mathbf{X}_n}(u)\,\mathrm{d}u
\nonumber\\
&&\quad=\int_{\mathbb{R}^{+}\times\mathbb{R}^{+}} \mathrm{e}^{-\psi_{\gamma+u}(z)}%
\tau_{n_{k}+1,Y_{k}}(\gamma+z+u)
\\&&\hspace*{51pt}{}
\times
\frac{%
\mathrm{e}^{-\psi_{0}(\gamma+u)}%
\prod_{k^{\prime}=1,k^{\prime}\neq k}^{\aleph(\mathbf{p}_{n})}%
\tau_{n_{k^{\prime}},Y_{k^{\prime}}}(\gamma+z+u)%
u^{n+q-1}}{\int_{\mathbb{R}^{+}}%
\mathrm{e}^{-\psi_{0}(\gamma+u)}%
\prod_{k^{\prime}=1}^{\aleph(\mathbf{p}_{n})}%
\tau_{n_{k^{\prime}},Y_{k^{\prime}}}(\gamma+u)u^{n+q-1}\,\mathrm{d}u}\,\mathrm{d}u\,\mathrm{d}z%
.
\nonumber
\end{eqnarray}
Taking the equality
$\psi_{0}(\gamma+z+u)=\psi_{\gamma+u}(z)+\psi_{0}(\gamma+u)$ and
transforming the upper integral with $(w,v)=(u+z,u)$,
(\ref{Prop1eqn11}) becomes
%
\begin{eqnarray}
\label{Prop1eqn12} &&\frac{\int_{\mathbb{R}^{+}}%
\mathrm{e}^{-\psi_{0}(\gamma+w)}%
\tau_{n_{k}+1,Y_{k}}(\gamma+w)%
\prod_{k^{\prime}=1,k^{\prime}\neq k}^{\aleph(\mathbf{p}_{n})}%
\tau_{n_{k^{\prime}},Y_{k^{\prime}}}(\gamma+w)%
{w^{n+q}}/({n+q})\,\mathrm{d}w}%
{\int_{\mathbb{R}^{+}}%
\mathrm{e}^{-\psi_{0}(\gamma+u)}%
\prod_{k^{\prime}=1}^{\aleph(\mathbf{p}_{n})}%
\tau_{n_{k^{\prime}},Y_{k^{\prime}}}(\gamma+u)u^{n+q-1}\,\mathrm{d}u}%
\nonumber\\[-8pt]\\[-8pt]
&&\quad=\frac{1}{n+q}\int_{\mathbb{R}^{+}}%
\frac{\tau_{n_{k}+1,Y_{k}}(\gamma+u)}%
{\tau_{n_{k},Y_{k}}(\gamma+u)}%
uf_{U_{n}\vert\mathbf{X}_{n}}(u)\,\mathrm{d}u .
\nonumber
\end{eqnarray}
We conclude that the predictive distribution for $X_{n+1}$ given
$\mathbf{X}_n= \{X_1,\ldots,X_n \}$ is given by
%
\begin{eqnarray}
\label{predX1} \mathbb{P} \{X_{n+1}\in \mathrm{d}x\vert\mathbf{X}_{n}
\} &=&%
\frac{1}{n+q}\int_{\mathbb{R}^{+}}u
\tau_{1,x}(\gamma+u) f_{U_n\vert\mathbf{X}_n}(u)\,\mathrm{d}u \times\eta(\mathrm{d}x)%
\nonumber\\[-8pt]\\[-8pt]
&&{}+ %
\sum_{k=1}^{\aleph(\mathbf{p}_n)}%
\frac{1}{n+q}\int_{\mathbb{R}^{+}}u\frac{\tau_{n_k+1,Y_k}(\gamma+u
)}{\tau_{n_k,Y_k}(\gamma+u)}f_{U_n\vert\mathbf{X}_n}(u)\,\mathrm{d}u%
\times\delta_{Y_k}(\mathrm{d}x)
\nonumber
.\quad
\end{eqnarray}
This representation is analogous to James \textit{et al.} \cite{JLP2008}. In fact, we
prefer an intuitive representation for this predictive distribution
for the next proof of the conditional case.
A representation of the
predictive distribution is given by
%
\begin{eqnarray}
\label{predX2} \mathbb{P} \{X_{n+1}\in \mathrm{d}x\vert\mathbf{X}_{n}
\} =%
\frac{\omega_{n,\aleph(\mathbf{p}_n)+1}(x)}%
{\phi(\mathbf{X}_n)}%
\eta(\mathrm{d}x)+%
\sum_{k=1}^{\aleph(\mathbf{p}_n)}%
\frac{\omega_{n,k}(Y_k)}%
{\phi(\mathbf{X}_n)}%
\delta_{Y_k}(\mathrm{d}x) ,
\end{eqnarray}
where
\begin{eqnarray*}
\omega_{n,\aleph(\mathbf{p}_n)+1}(x)&=&\int_{\mathbb{R}^{+}}u\tau_{1,x}(
\gamma+u)f_{U_n\vert\mathbf{X}_n}(u)\,\mathrm{d}u ,
\\
\omega_{n,k}(Y_k)&=&\int_{\mathbb{R}^{+}}u
\frac{\tau_{n_k+1,Y_k}(\gamma
+u)}{\tau_{n_k,Y_k}(\gamma+u)}f_{U_n\vert\mathbf{X}_n}(u)\,\mathrm{d}u ,
\\
\phi(\mathbf{X}_n)&=&\int_{\mathbb{X}}
\omega_{n,\aleph(\mathbf
{p}_n)+1}(x)\eta(\mathrm{d}x)+\sum_{k=1}^{\aleph(\mathbf{p}_n)}
\omega_{n,k}(Y_k)=n+q .%
\end{eqnarray*}
Clearly, we could write (\ref{predX1}) or (\ref{predX2}) as an
expectation with respect to the distribution $U_n$ given
$\mathbf{X}_n$,
$\mathbb{P} \{X_{n+1}\in \mathrm{d}x\vert\mathbf{X}_{n} \}= E [\mathbb
{P} \{X_{n+1}\in \mathrm{d}x\vert\mathbf{X}_{n},U_n \}\vert\mathbf
{X}_{n} ]$. However, the term $\mathbb{P} \{X_{n+1}\in \mathrm{d}x\vert
\mathbf{X}_{n},U_n \}$ might not be a proper distribution. We
suggest the following expectation instead with respect to the
distribution $\widetilde{U}_n$ given $\mathbf{X}_n$,
$\mathbb{P} \{X_{n+1}\in \mathrm{d}x\vert\mathbf{X}_{n} \}= E [\mathbb
{P} \{X_{n+1}\in \mathrm{d}x\vert\mathbf{X}_{n},\widetilde{U}_n \}\vert
\mathbf{X}_{n} ]$. So conditional on $\widetilde U_n$, the
predictive distribution is given by
\begin{eqnarray*}
\mathbb{P} \{X_{n+1}\in \mathrm{d}x\vert\mathbf{X}_{n},
\widetilde{U}_n \} &=&%
\frac{\omega_{n,\aleph(\mathbf{p}_n)+1}(\widetilde{U}_n,x)}{\phi
(\widetilde{U}_n,\mathbf{X}_n)}%
\eta(\mathrm{d}x)+%
\sum_{k=1}^{\aleph(\mathbf{p}_n)}%
\frac{\omega_{n,k}(\widetilde{U}_n,Y_k)}{\phi(\widetilde{U}_n,\mathbf
{X}_n)}%
\delta_{Y_k}(\mathrm{d}x) ,
\\
\phi(u,\mathbf{X}_n) &=& \int_{\mathbb{X}}
\omega_{n,\aleph(\mathbf{p}_n)+1}(u,x)\eta(\mathrm{d}x)+\sum_{k=1}^{\aleph(\mathbf{p}_n)}
\omega_{n,k}(u,Y_k) ,
\\
\omega_{n,\aleph(\mathbf{p}_n)+1}(u,x)&=&u\tau_{1,x}(\gamma+u) ,
\\
\omega_{n,k}(u,Y_k)&=&u \frac{\tau_{n_k+1,Y_k}(\gamma+u)}{\tau_{n_k,Y_k}(\gamma+u)}
,\qquad
k=1,\ldots,\aleph(\mathbf{p}_n) \quad \mbox{and}
\\
f_{\widetilde{U}_n\vert\mathbf{X}_n}(u)&=&%
\frac{\phi(u,\mathbf{X}_n)}{E[\phi(U_n,\mathbf{X}_n)\vert\mathbf
{X}_n]}f_{U_n\vert\mathbf{X}_n}(u) .
\end{eqnarray*}
Lastly, the following equality can be achieved by inspection of
(\ref{predX1}) and (\ref{predX2}),\linebreak[4]
$\phi(\mathbf{X}_n)= E[\phi(U_n\vert\mathbf{X}_n)\vert\mathbf
{X}_n]=n+q$
and the following equalities can be obtained according to the
definitions,
$\omega_{n,\aleph(\mathbf{p}_n)+1}(x)=E[\omega_{n,\aleph(\mathbf
{p}_n)+1}(U_n,x)\vert\mathbf{X}_n]$ and
$\omega_{n,k}(Y_k)=E[\omega_{n,k}(U_n,Y_k)\vert\mathbf{X}_n]$
for $k=1,\ldots,\aleph(\mathbf{p}_n)$.

\subsection{\texorpdfstring{Proof of Proposition \protect\ref{prop2marginaldistribution}}
{Proof of Proposition 3.2}}\label{prop2proof}
The joint distribution of
$ \{\mathbf{Y}_{\aleph(\mathbf{p}_n)},\mathbf{p}_n,U_n \}$
could be simply derived from
(\ref{Thm2eqn4}) of Theorem \ref{thm2conditionalposteriorG}, which
is proportional to
\begin{eqnarray*}
\mathrm{e}^{-\psi_{0}(\gamma+u)}%
\frac{u^{n+q-1}}{\Gamma(n+q)} \prod
_{k=1}^{\aleph(\mathbf{p}_{n})}%
\eta(\mathrm{d}Y_{k})
\tau_{n_{k},Y_k}(\gamma+u)%
\,\mathrm{d}u,
\end{eqnarray*}
where $\tau_{m,z}(a)$ is defined in (\ref{tau}).
Distributional results around these random variables $ \{\mathbf
{Y}_{\aleph(\mathbf{p}_n)},\mathbf{p}_n,U_n \}$
are achieved immediately with some simple algebra.

\subsection{\texorpdfstring{Proof of equation (\protect\ref{partitiondistributiongivent})}
{Proof of equation (6.1)}}\label{partitiondistributiongiventproof}
Starting from the marginal distribution of $\mathbf{X}_n$ which is
given by
%
\begin{eqnarray}
\label{PGTeqn1} \int_{\mathcal{M}_\mathbb{X}}%
\prod
_{i=1}^{n}%
\frac{\mu(\mathrm{d}x_{i})}{\mu(\mathbb{X})}%
\mathcal{P}_{\widetilde\mu_t}(\mathrm{d}\mu) =%
\int_{\mathcal{M}_\mathbb{X}}%
\prod_{i=1}^{n}%
\mu(\mathrm{d}x_{i})%
\frac{h(\mu(\mathbb{X}))\mu(\mathbb{X})^{-n}}{E[h(\mu(\mathbb{X}))]} \mathcal{P}_{\widetilde\mu}(\mathrm{d}
\mu) .%
\end{eqnarray}
One could apply the Fubini theorem following from an application of
Lemma 2.2 of James \cite{James2002}, page 8, on (\ref{PGTeqn1}) to yield the
joint distribution of $ \{\mathbf{Y}_{\aleph(\mathbf{p}_n)},\mathbf
{p}_n \}$,
%
\begin{eqnarray}
\label{PGTeqn2} \frac{1}{E[h(\mu(\mathbb{X}))]} 
\int
_{\mathcal{M}_\mathbb{X}\times(\mathbb{R}^+)^{\aleph(\mathbf{p}_n)}} 
\frac{h(\mu(\mathbb{X})+\sum_{k=1}^{\aleph(\mathbf{p}_n)}s_k)}%
{(\mu(
\mathbb{X})+\sum_{k=1}^{\aleph(\mathbf{p}_n)}s_k
)^{n}} \mathcal{P}_{\widetilde\mu}(\mathrm{d}\mu) \prod
_{k=1}^{\aleph(\mathbf{p}_n)}%
\eta(\mathrm{d}Y_k)
s_k^{n_k} \rho(\mathrm{d}s_k ) .\qquad
\end{eqnarray}
Then, in (\ref{PGTeqn2}), marginalizing $\mathbf{Y}_{\aleph(\mathbf
{p}_n)}$ and
replacing the distribution of the total $\widetilde T=\widetilde\mu
(\mathbb{X})$ by
$f_{\widetilde T}(t)$, the distribution of $\mathbf{p}_n$ is given by
the following without the proportional constant,
%
\begin{eqnarray}
\label{PGTeqn3} 
\int_{(\mathbb{R}^+)^{\aleph(\mathbf{p}_n)}\times\mathbb{R}^+}%
\frac{h(t+\sum_{k=1}^{\aleph(\mathbf{p}_n)}s_k)}%
{(t+\sum_{k=1}^{\aleph(\mathbf{p}_n)}s_k
)^{n}} \frac{1}{E[h(\widetilde T)]} f_{\widetilde T}(t) \prod
_{k=1}^{\aleph(\mathbf{p}_n)}%
s_k^{n_k}
\rho(\mathrm{d}s_k )\, \mathrm{d}t .
\end{eqnarray}
Considering the transformation $t^*=t+\sum_{k=1}^{\aleph(\mathbf
{p}_n)}s_k$ and (\ref{PGTeqn3}) becomes
%
\begin{eqnarray}
\label{PGTeqn4} 
\int_{(\mathbb{R}^+)^{\aleph(\mathbf{p}_n)}\times\mathbb{R}^+}%
\mathbb{I}_{\{t^*-\sum_{k=1}^{\aleph(\mathbf{p}_n)}s_k>0\}} \frac{h(t^*)}%
{(t^*)^{n}}%
\frac{f_{\widetilde T}(t^*-\sum_{k=1}^{\aleph(\mathbf
{p}_n)}s_k)}{E[h(\widetilde T)]} \prod
_{k=1}^{\aleph(\mathbf{p}_n)}%
s_k^{n_k}
\rho(\mathrm{d}s_k )\, \mathrm{d}t^* .
\end{eqnarray}
Rearranging terms in (\ref{PGTeqn4}) yields
\begin{eqnarray*}
\int_{(\mathbb{R}^+)^{\aleph(\mathbf{p}_n)}\times\mathbb{R}^+}%
\mathbb{I}_{\{t^*-\sum_{k=1}^{\aleph(\mathbf{p}_n)}s_k>0\}}
\frac{f_{\widetilde T}(t^*-\sum_{k=1}^{\aleph(\mathbf{p}_n)}s_k)}%
{(t^*)^{n}f_{\widetilde T}
(t^*)}%
\prod_{k=1}^{\aleph(\mathbf{p}_n)}s_k^{n_k}
\rho(\mathrm{d}s_k)%
\frac{h(t^*)f_{\widetilde T}(t^*)}{E[h(\widetilde T)]}%
\,\mathrm{d}t^*.
\end{eqnarray*}
Thus, the proof is complete.
\end{appendix}

\section*{Acknowledgements}

The author thanks Edward Cripps for the tremendous help in revising this
article and Robin Milne for proofreading this manuscript and
giving constructive comments on the presentation.
Their invaluable comments and supports to the realization of this work
are deeply appreciated.
The author also thanks the associate editor and two
referees for their helpful suggestions concerning the presentation
of this article and recommending various improvements in exposition.
This research was partially supported by Australian Actuarial Research
Grant \#9201101267.



\printhistory

\end{document}